\documentclass[11 pt]{article}
\usepackage{amssymb}
\usepackage{latexsym}
\usepackage{amsmath}
\usepackage{epsfig}
\usepackage{verbatim}
\usepackage{makeidx}

\newcommand{\IP}{\mathbb P}
\newcommand{\IQ}{\mathbb Q}
\newcommand{\IE}{\mathbb E}
\newcommand{\IR}{\mathbb R}

\newcommand{\IZ}{\mathbb Z}
\newcommand{\IN}{\mathbb N}

\newcommand{\la}{\langle}
\newcommand{\ra}{\rangle}
\newtheorem{thm}{Theorem}[section]
\newtheorem{propn}[thm]{Proposition}
\newtheorem{lemma}[thm]{Lemma}

\newtheorem{remark}[thm]{Remark}
\newtheorem{cor}[thm]{Corollary}

\newenvironment{proof}{\noindent
PROOF.}{\hfill$\mathbf{\Box}$ \vspace{.5\baselineskip}}
\DeclareMathOperator{\sgn}{sgn}

\begin{document}

\title{On pathwise uniqueness for stochastic heat equations
with non-Lipschitz coefficients}

\author{Leonid Mytnik$\hbox{}^{1}$\\
Faculty of Industrial Engineering and Management\\
Technion Israel Institute of Technology\\Haifa 32000\\ Israel\\
leonid@ie.technion.ac.il\\~\\
Edwin Perkins$\hbox{}^{2}$\\
Mathematics Department\\
The University of British Columbia\\
1984 Mathematics Road\\
Vancouver, British Columbia, Canada V6T 1Z2\\
perkins@math.ubc.ca\\~\\
Anja Sturm$\hbox{}^{3}$\\
Department of Mathematical Sciences\\
University of Delaware\\
501 Ewing Hall\\
Newark, Delaware, 19716-2553, \\
U.S.A.\\
sturm@math.udel.edu}

\maketitle
\begin{abstract}
We consider the existence and pathwise uniqueness of the
stochastic heat equation with a multiplicative colored
noise term on $\IR^{d}$ for $d \geq 1.$ We focus on the case of
non-Lipschitz noise coefficients and singular spatial noise correlations.
In the course of the proof a new result on H\"older continuity
of the solutions near zero is established.
\end{abstract}
        \vspace{\fill}

          \thispagestyle{empty}

\par
         (1) Research supported in part by the Israel Science Foundation
(grant No. 116/01 - 10.0).\\
\par        (2) Research supported in part by an NSERC Research grant.\\

\par (3) Research supported in part by an NSERC Research grant and by DFG Priority Program 1033.

\newpage

\section{Introduction}
\label{section:intro}

\noindent

This work is motivated by the following question: Does pathwise
uniqueness hold in the parabolic stochastic  pde
\begin{equation}
\label{sbm}
\frac{\partial}{\partial t}
u(t,x) = \frac{1}{2} \Delta u(t,x) dt + \sqrt{u(t,x)} \dot{W}(x,t)?
\end{equation}
Here $\Delta$ denotes the Laplacian and $\dot{W}$ is space-time white
noise on $\IR_+\times \IR$. It is known that uniqueness in law holds for
solutions to (\ref{sbm}) in the appropriate space of continuous functions
and such solutions are the density for one-dimensional super-Brownian
motion (see, e.g., Section III.4 of \cite{P02}).  One motivation for
studying pathwise uniqueness is the hope that such an approach would be
more robust and establish uniqueness for closely related equations in
which $\sqrt{u(t,x)}$ could be replaced by $\sqrt
{\gamma(u(t,x))u(t,x)}$. Such models arise as scaling limits of
critical branching particle systems in which the branching rate at
$(t,x)$ is given by $\gamma(u(t,x))$.  The method used to establish uniqueness
in law for solutions of (\ref{sbm}) is duality. This approach has the advantage
of giving a rich toolkit for the study of solutions to (\ref{sbm}) but the
disadvantage of being highly non-robust, although one of us was able to extend
this method to powers of $u(t,x)$ between $1/2$ and $1$ (see \cite{M98}).  

The difficulty in proving pathwise uniqueness in (\ref{sbm}) arises from the
fact that
$\sqrt{u}$ is non-Lipschitz.  The above equation does have the advantage 
of having a diagonal form--that is, when viewed as a continuum-dimensional
stochastic differential equation there are no off-diagonal terms in the
noise part of the equation and the diffusion coefficient for the $x$ coordinate is a
function of that coordinate alone.  For finite-dimensional sde's this was the setting for
Yamada and Watanabe's extension (\cite{YW71}) of Ito's pathwise uniqueness results to
H\"older continuous coefficients, and so an optimist may hope this approach can carry over
to our infinite dimensional setting.  As we will be using their
conditions later, let us recall the Yamada-Watanabe result. Let $\rho$ be a strictly
increasing function on
$\IR_{+}$ such that
\begin{equation}
\label{rhocond}
\int_{0+} \rho^{-2}(x) dx = \infty.
\end{equation}
Now assume that $\sigma:\IR\to\IR$ is such that for all $x,y \in \IR,$
\begin{equation}
\label{cond:Yamada}
|\sigma(x) - \sigma(y)| \leq \rho(|x-y|).
\end{equation} 
Then pathwise uniqueness holds for solutions of the one-dimensional sde
\begin{equation}
\label{ywsde}
X(t)=X(0)+\int_0^t\sigma(X(t))dB(t),
\end{equation}
where $B$ is a standard Brownian motion.
The square root function clearly satisfies the above hypotheses but the
infinite dimensional setting has stymied attempts to carry the
methodology over.  Yamada and Watanabe's proof has been simplified (see
e.g., Theorem IX.3.5 of \cite{RY91}) by the notion of the local time of
a semimartingale and the fact that $u(t,x)$ will not be a semimartingale
in $t$ for $x$ fixed (it will only be H\"older continuous of index $1/4$)
would seem to be a serious obstacle in directly applying these
methods.  

We will not resolve the uniqueness question posed above, but will succeed
in extending the above ideas to stochastic heat equations of the
form
\begin{equation}
\label{heatcol}
\frac{\partial}{\partial t}
u(t,x) = \frac{1}{2} \Delta u(t,x) dt + \sigma(u(t,x)) \dot{W}(x,t).
\end{equation}
for colored noises other than white, and appropriate H\"older continuous,
but not necessarily Lipschitz continuous,
$\sigma$.  Here, $u$ is a random function on
$\IR_{+} \times \IR^{d}$ and
we sometimes write $u_{t}$ for $u(t, \cdot).$
The coefficient $\sigma$ is a real-valued continuous function on
$\IR$. It is assumed throughout this work to satisfy
the following global growth condition: For all $u \in \IR$
there exists a constant $c_{\ref{growthcond}}$ such that
\begin{equation}
\index{linear growth condition}
\label{growthcond}
|\sigma(u)| \leq c_{\ref{growthcond}} (1 + |u|).
\end{equation}
Here and elsewhere $c_i$ and $c_{i.j}$ will denote fixed positive constants,
while $C$ will denote a positive constant which may change from line to line.
\noindent
The noises $W$ considered here are Gaussian martingale measures
on $\IR_{+} \times \IR^{d}$ in the sense of Walsh \cite{jW}. $W$ is defined
on a filtered probability space $(\Omega, {\cal F}, {\cal F}_{t}, \IP)$ and
$W_t(\phi)=\int_0^t\int_{\IR^d} \phi(s,x) W(dxds)$ is an ${\cal
F}_t$-martingale for $\phi\in C_{c}^{\infty}(\IR_{+} \times \IR^{d}),$ the
space of compactly supported, infinitely differentiable functions  on
$\IR_+\times\IR^d$. If $W(\phi)=W_\infty(\phi)$, $W$ can be characterized by
its covariance functional
\begin{equation}
\label{knoise}
J_{k}(\phi,\psi):= \IE\left[ W(\phi) W(\psi)\right]
= \int_{0}^{\infty}\int_{\IR^{d}} \int_{\IR^{d}}
\phi(s,x) k(x,y) \psi(s,y) dx dy ds,
\end{equation}
for $\phi,\psi
\in C_{c}^{\infty}(\IR_{+} \times \IR^{d})$.  We call the function
\label{correlationkernel} $k:\IR^{2d} \rightarrow \IR$
the correlation kernel of $W.$ Some sufficient
conditions for the existence of a martingale measure
$W$ corresponding to $k$ are that $J_{k}$ is
symmetric, positive definite and continuous. Thus, necessarily,
$k(x,y)=k(y,x)$ for all $x,y \in \IR^{d}.$ Continuity on
$C_{c}^{\infty}$ is implied, for example, if $k$ is integrable
on compact sets. We also note that a general class of martingale
measures, spatially homogeneous
noises, can be described by (\ref{knoise}) where
$k(x,y)=\tilde{k}(x-y)$. 

If $\sigma(u)=u$ then equation (\ref{heatcol}) arises as the diffusion limit
of super-Brownian motion in $\IR^d$ where the offspring law depends on a random
environment, whose spatial correlation is described by $k.$ For $k$ bounded, this was
proven in \cite{aS03}. More general coefficients $\sigma$ may be thought of as reflecting
an additional dependence of the offspring law on the local particle density.
 
If $k$ is bounded,   Viot \cite{mV76} proved pathwise uniqueness for solutions to
(\ref{heatcol}) on bounded domains of $\IR^{d}$ for
$\sigma(u)= \sqrt{u (1-u)}_{+},$ where the subscript indicates that the
positive part of the function is taken.  We will extend this result to
our setting for solutions of (\ref{heatcol}) on $\IR^d$ with bounded $k$
in Theorem~\ref{thm:uniquekbounded} below.  Note that white noise will
correspond to the case where we set $\tilde k$ equal to the generalized
function $\delta_0$ in the above.  Our main result
(Thm~\ref{thm:uniquekunbounded} below) will interpolate between these
settings and establish pathwise uniqueness for colored noises for which
the correlation is bounded by a Riesz kernel,
\begin{equation}
\label{def:singulartildek}
|k(x,y)| \leq c_{\ref{def:singulartildek}} [|x-y|^{-\alpha}+1]\hbox{ for
all
$x,y\in\IR^d$ and appropriate }
\alpha>0.
\end{equation}

In order to formulate a condition on the singularity of ${k}$ and
relate our conditions to those in the literature, we define the spectral
measure, $\mu$, of a spatially homogeneous covariance kernel $\tilde k$:
\begin{equation}
\label{eq:specmeas}
\int_{\IR^{d}} \tilde{k}(x) \phi(x) dx = \int_{\IR^{d}} {\cal F} \phi(\xi) \mu(d\xi)
\end{equation}
for any rapidly decreasing test function $\phi$ where
${\cal F} \phi(\xi) = \int_{\IR^{d}} \exp( -2 i \pi \xi \cdot x) \phi(x) dx$
is the Fourier transform.
Later on we will assume $\mu$ to be a tempered measure fulfilling for
some $\eta \in [0,1],$
\begin{equation}
\label{eq:colorcond}
\int_{\IR^{d}} \frac{\mu(d\xi)}{(1+|\xi|^{2})^{\eta}} < \infty.
\end{equation}
To relate (\ref{def:singulartildek}) with condition (\ref{eq:colorcond}) 
used in the literature, we introduce:
\begin{itemize}
\item[$(A)_\eta$:] ($\eta>0$)
$W$ is a Gaussian noise with correlation kernel $|k(x,y)| \leq
c_{\ref{eq:colorcond}}
\tilde{k}(x-y), x,y \in \IR^{d}$ for some symmetric, locally bounded and
positive definite kernel $\tilde k$ whose spectral measure satisfies
(\ref{eq:colorcond}).
\item[$(A)_0$:] $W$ is a Gaussian noise and its correlation
kernel $k$ is bounded. \\
\end{itemize}

\begin{remark}\label {rem:sm}
Note that (\ref{def:singulartildek}) implies $(A)_\eta$ for $\alpha \in (0, 2\eta \wedge d):$
Here, $\tilde{k}(x) =  |x|^{-\alpha}+1$ and the spectral measure is of the
form
$\mu(d\xi) = c_{\ref{rem:sm}} [|\xi|^{\alpha-d}d\xi+\delta_0(d\xi)]$. Hence,
condition (\ref{eq:colorcond}) is satisfied if and only if $\alpha \in
(0, 2\eta
\wedge d)$ (see Chap. V Lemma 2(a) of \cite{eS67}).  Note also that the positive definite
spatially homogeneous kernels $k_\alpha(x,y)=|x-y|^{-\alpha}$ give a natural family of 
kernels for which our results will hold.
\end{remark}

In order to make sense of the formal equation (\ref{heatcol}) we use
the variation of constants form of solutions:
Denote by $p$ be the $d$-dimensional heat kernel
\begin{equation}
  \label{eq:heatkernel}
  p_{t}(x) = \frac{1}{(2 \pi t)^{\frac{d}{2}}}
\exp(-\frac{|x|^{2}}{2t}).
\end{equation}

A stochastic process $u: \Omega \times \IR_{+} \times \IR^{d}
\rightarrow \IR$, which is jointly measurable and ${\cal F}_{t}$-adapted,
is said to be a solution to the stochastic heat equation (\ref{heatcol})
in the variation of constants sense with respect to the martingale
measure $W,$ defined on $\Omega,$
and initial condition $u_{0}$, if for each $t\ge 0$, a.s. for almost
all
$x
\in
\IR^{d}$
\begin{eqnarray}
\label{heatecolint}
u(t,x) &=& \int_{\IR^{d}} p_{t}(x-y) u_{0}(y) dy
+ \int_{0}^{t}  \int_{\IR^{d}} p_{t-s}(x-y) \sigma(u(s,y)) W(dy ds).
\end{eqnarray}
\noindent

Solutions to (\ref{heatecolint}) have been well studied in the case where
$\sigma$ is Lipschitz continuous in $u.$ A sufficient condition
for strong existence and uniqueness
of solutions is given by $(A)_{\eta}$ for
$\eta \leq 1,$ see Dalang
\cite{rD99} (see also Theorem \ref{thm:Dalang} in the Appendix) and Peszat and Zabczyk
\cite{PZ00}. H\"older continuity of the sample paths was
established by Sanz-Sol\'e
and Sarr\`a \cite{SS02} if $\eta <1$ (cf. 
Lemma \ref{thm:SanzSoleSarra} in the
Appendix).

To state the main results we introduce some notation, which will
be used throughout this work:  We write $C(\IR^{d})$ for the space of continuous
functions on $\IR^{d}.$ A superscript $k,$ respectively $\infty$, indicates that
functions are in addition $k$ times, respectively infinitely often, continuously
differentiable. A subscript $b,$ respectively $c,$ indicates that they are also
bounded, respectively have compact support. We also define
\begin{equation*}
||f||_{\lambda,\infty}:=\sup_{x \in \IR^{d}} |f(x)| e^{-\lambda |x|},
\end{equation*}
set $C_{tem}:=\{f \in C(\IR^{d}) , ||f||_{\lambda,\infty}
< \infty  \text{ for any } \lambda >0\}$ and endow it with the topology induced
by the norms $||\cdot ||_{\lambda,\infty}$ for $\lambda>0.$ That is,
$f_n\to f$ in $C_{tem}$ iff $\lim_{n\to\infty}\Vert
f-f_n\Vert_{\lambda,\infty}=0$ for all $\lambda>0$.  For
$I\subset \IR_+$, let
$C(I,E)$ be the space of all continuous functions on $I$ taking values in
a topological space $E$, endowed with the topology of uniform convergence
on compact subsets of $I$.  A stochastically weak solution to
(\ref{heatecolint}) is a solution on some filtered space with respect to some
noise $W$, i.e., the noise and space is not specified in advance.

With this notation we can state the following standard existence
result whose proof is outlined in the Appendix:
\begin{thm}
\label{thm:nonLipschitzexistence}
Let $u_{0} \in C_{tem},$ and  let $\sigma$ be a
continuous function satisfying the growth bound
(\ref{growthcond}). Assume that (\ref{def:singulartildek}) holds for
some $\alpha \in (0,2\wedge d).$  Then there exists a stochastically
weak solution to (\ref{heatecolint}) with sample paths a.s. in
$C(\IR_{+}, C_{tem})$.
\end{thm}

\begin{remark}
\label{remark:nonLipschitzex2}(a) The proof in fact only requires that
$(A)_\eta$ hold for some $\eta\in[0,1)$, a condition which follows from the
above bound on $k$ by Remark~\ref{rem:sm}.

\medskip 
\noindent(b) In the case where the correlation kernel is bounded,
existence has been shown for more general initial conditions and solution spaces
in Sturm \cite{aS03}: Define $L^{p}_{\lambda}(\IR^{d}):=L^{p}(\IR^{d},e^{-\lambda |x|} dx)$
and denote the associated norm by $|| \cdot||_{\lambda,p}.$ Then if
$\IE ( ||u_{0}||^{p}_{\lambda,p} )< \infty,$
for some $p > 2$ and $\lambda>0 $ there exists a stochastically weak solution
$u \in C(\IR_{+},L^{p}_{\lambda}(\IR^{d})),$
to (\ref{heatecolint}) which satisfies
\begin{equation}
\label{uniformbound}
 \IE \left( \sup_{0 \leq t \leq T}
||u(t,\cdot)||^{p}_{\lambda,p} \right) <\infty \hbox{  for any $T>0$.}
\end{equation}

\end{remark}

We say pathwise uniqueness holds for solutions of (\ref{heatecolint}) in
$C(\IR_+,C_{tem})$ if for every $u_0\in C_{tem}$, any two solutions to
(\ref{heatecolint}) with sample paths a.s. in $C(\IR_+,C_{tem})$ must be equal
with probability
$1$.  For Lipschitz continuous $\sigma$, it is easy to modify
Theorem~13 of \cite{rD99} and Theorem~2.1 of
\cite{SS02} to get pathwise uniqueness and H\"older continuity of solutions for
$\alpha<2\wedge d$.  Also, Theorem~11 and Remark~12 of
\cite{rD99} show that function-valued solutions will not exist for
$\alpha>2\wedge d$.   Here then is our --it holds in any spatial
dimension
$d$:

\begin{thm}
\label{thm:uniquekunbounded}
Assume that for some $\alpha\in (0,1)$, $\sigma:\IR\to\IR$ satisfies 
(\ref{growthcond}), is H\"older continuous of index $\gamma$ for some
$\gamma\in
(\frac{1+\alpha}{2} , 1]$, and 
\[|k(x,y)|\leq c_{\ref{thm:uniquekunbounded}}[|x-y|^{-\alpha}+1]
\hbox{ for all } x,y \in \IR^{d}.\]
Then pathwise uniqueness holds for 
solutions of  (\ref{heatecolint}) in $C(\IR_{+},C_{tem})$.
\end{thm}

\begin{remark} The H\"older condition on $\sigma$ may be weakened to 
to the local H\"older condition: \\
For any $K>0$ there exists $L=L(K)$ such that 
\[|\sigma(u)-\sigma(v)|\le L(|u-v|^\gamma+|u-v|)\ \ \forall u,v: \;|u|,|v|\leq K,\]
where $\gamma$ is as Theorem~\ref{thm:uniquekunbounded}.
The required modifications in the proof are elementary.  
\end{remark}
In the above result there is a trade-off between the H\"older continuity of
$\sigma$ and the singularity of the covariance kernel of the noise.  For $d=1$, letting
$\alpha\to1-$ and renormalising will give white noise.  More specifically, if $\tilde
k_\alpha(x-y)={1-\alpha\over 2}|x-y|^{-\alpha}$, then  for
$\phi,\psi\in C_c^\infty(\IR_+\times\IR^d)$, $\lim_{\alpha\to1-}J_{\tilde
k_\alpha}(\phi,\psi)=\int_0^\infty \int \phi(s,x)\psi(s,x)dxds$.  
The H\"older condition in Theorem~\ref{thm:uniquekunbounded} 
approaches Lipschitz continuity. (As $\tilde k$ should be locally integrable we cannot
expect to take $\alpha=1$.) Hence, although the result does not say anything about white
noise itself, it at least coincides with the known Lipschitz conditions which imply
pathwise uniqueness in the limit as $\alpha$ approaches $1$.  The same cannot be said for
higher dimensions.  Here, the aforementioned results of Dalang, and Sanz-Sol\'e and
Sarr\'a show that for $\alpha< 2$ we will
have pathwise unique continuous solutions when the coefficients are Lipschitz
continuous.  Unfortunately, our hypotheses become vacuous in the above
uniqueness theorem when
$\alpha$ exceeds $1$ and so we believe our condition on the H\"older index in
Theorem~\ref{thm:uniquekunbounded} is non-optimal in dimensions greater
than $1$.  At the other end of the scale we see that as
$\alpha$ approaches $0$, the required H\"older exponent approaches
$1/2$, the critical power in the one-dimensional results of Yamada and
Watanabe.  In fact, if the covariance kernel is bounded we can
weaken the H\"older condition on 
$\sigma$ to precisely the Yamada-Watanabe condition
(\ref{rhocond},\ref{cond:Yamada}) introduced above. Again the result holds in
any spatial dimension.

\begin{thm}
\label{thm:uniquekbounded}
Assume that $(A)_{0}$ holds and that $\sigma:\IR\to\IR$ satisfies
(\ref{growthcond}) and  (\ref{cond:Yamada}).  Then pathwise uniqueness holds
for solutions of (\ref{heatecolint}) in $C(\IR_{+},C_{tem})$.
\end{thm}

\begin{remark} (a) The conclusions of Theorems \ref{thm:nonLipschitzexistence},
\ref{thm:uniquekunbounded} and \ref{thm:uniquekbounded} remain valid if we
allow for an additional drift term in the heat equation. More precisely, we can
add a term of the form
$\int _0^t\int p_{t-s}(x-y)f(u(s,y))dyds$ to the right hand side
of (\ref{heatecolint}), where $f$ satisfies the growth bound
(\ref{growthcond}), is continuous in the existence theorem,
Theorem~\ref{thm:nonLipschitzexistence}, and is Lipschitz continuous for the
uniqueness results, Theorems~\ref{thm:uniquekunbounded} and
\ref{thm:uniquekbounded}. The additional arguments are standard.
\medskip

\noindent(b) The pathwise uniqueness conclusions of Theorems~~\ref{thm:uniquekunbounded} and
\ref{thm:uniquekbounded}, and weak existence given by
Theorem~\ref{thm:nonLipschitzexistence} imply the existence of a strong solution
to (\ref {heatecolint}), that is a solution which is adapted with respect to the
canonical filtration of the noise $W$. The proof follows just as in the
classical sde argument of Yamada and Watanabe (see, e.g. Theorem IX.1.7 of
\cite{RY91}).
\medskip

\noindent(c) Theorem~\ref{thm:uniquekbounded} holds true if we consider
solutions with paths in $C(\IR_{+}, L^{p}_{\lambda}(\IR^{d}))$ as was done in
Viot's work \cite{mV76}. In fact, the arguments given in sections \ref{section:auxiliary}
and \ref{section:proofkbounded} remain the same in this case. The only difference
is that a bit more care has to be taken to justify some of the convergences 
as the solutions are not necessarily continuous. But this can be done in a 
straightforward way. 
\end{remark}

The proof of our pathwise uniqueness theorems will require some moment bounds
for {\it arbitrary} continuous 
$C_{tem}$-valued solutions to the equation~(\ref{heatecolint}). Let
$S_t\phi(x)=\int p_t(y-x)\phi(y)\,dy$.  The following result will be proved in
the Appendix. 
\begin{propn}
\label{prop}
Let $u_{0} \in C_{tem},$ and  let $\sigma$ be a
continuous function satisfying the growth bound
(\ref{growthcond}). Assume that (\ref{def:singulartildek}) holds for
some $\alpha \in (0,2\wedge d).$  Then any solution
$u\in C(\IR_{+}, C_{tem})$ to (\ref{heatecolint})
has  the following properties.
\begin{itemize}
\item[{\bf (a)}]
For any $T,\lambda >0$ and $p \in (0,\infty),$
\begin{equation}
\label{eq:uniformmoments2a}
\IE\Big( \sup_{0 \leq t \leq T} \sup_{x \in \IR^{d}}
|u(t,x)|^{p} e^{-\lambda |x|}   \Big)
< \infty.
\end{equation}
\item[{\bf (b)}]
For any   $\xi \in (0, 1-\alpha/2)$ the process $u(\cdot,\cdot)$ is a.s.
uniformly H\"older continuous on compacts in $(0,\infty)\times \IR^{d}$, and 
the process $Z(t,x)\equiv u(t,x)-S_tu_0(x)$ is uniformly H\"older continuous on
compacts in $[0,\infty)\times \IR^d$, both with H\"older coefficients
$\frac{\xi}{2}$ in time and
$\xi$ in space.

Moreover, for any $T, R>0,$ and
$0\le t,t' \leq T, x,x' \in \IR^{d}$ such that $|x-x'|<R$ as well as $p \in [2,
\infty)$ and $\xi \in (0, 1-\alpha/2),$ there exists a constant
$c_{\ref{eq:Kolmogorov}}=c_{\ref{eq:Kolmogorov}}(T,p,\lambda, R,\xi)$ such that
\begin{equation}
\label{eq:Kolmogorov}
\IE\left( |Z(t,x) - Z(t',x') |^{p}e^{-\lambda|x|}\right)
\leq  c_{\ref{eq:Kolmogorov}}\left(|t-t'|^{\frac{\xi}{2} p}
+ |x-x'|^{\xi p} \right).
\end{equation}
\end{itemize}
\end{propn}

\begin{remark} The proof of the above will only require $(A)_\eta$ for some
$\eta\in[0,1)$, a condition which is implied by the hypotheses above (see
Remark~\ref{rem:sm}).  In this case we should take $\xi\in (0,1-\eta)$ in (b)
as is done in the proof in the Appendix.
\end{remark}

It is straightforward to show that under the hypotheses of Theorem
\ref{thm:nonLipschitzexistence}, solutions to (\ref{heatecolint}) with
continuous $C_{tem}$-valued paths are also solutions to the heat equation in its
distributional form for suitable test functions
$\Phi$. More specifically, for $\Phi \in C_{c}^{\infty}(\IR^{d}):$
\begin{eqnarray}
\label{weakheatb}
\int_{\IR^{d}} u(t,x) \Phi(x)dx
&=& \int_{\IR^{d}} u_{0}(x)  \Phi(x)dx
+\int_{0}^{t}\int_{\IR^{d}} u(s,x) \frac{1}{2} \Delta\Phi(x) dx ds
\\
\nonumber
& &+\int_{0}^{t}\int_{\IR^{d}}\sigma(u(s,x))\Phi(x) W(dx ds)
\quad\forall t\ge 0 \quad a.s.
\end{eqnarray}
In fact, given an appropriate class of test functions,
the two notions of solution (\ref{heatecolint}) and (\ref{weakheatb})
are equivalent. In our case, $\{ \Phi \in C^{\infty}(\IR^{d}): \Phi(x)
\leq C e^{-\lambda |x|} \text{ for some } C>0 \text{ and all } x \in \IR^{d} \}$ is a
suitable class of test functions.  For the
details of the proof we refer to Sturm \cite{aS02} Proposition 3.2.3. There, 
the setting is a bit different as it works in the setting 
of Remark \ref{remark:nonLipschitzex2} with bounded $k$. However, the
arguments do not change for the case of $k$ unbounded as long as the stochastic
integral in (\ref{weakheatb}) is well defined, which can easily be checked.

We now briefly outline the proof of our main result
(Theorem~\ref{thm:uniquekunbounded}) and the contents of the paper. To emulate
Yamada and Watanabe, consider a pair of solutions, $u^1$ and $u^2$, to
(\ref{heatecolint}), set $\tilde u=u^1-u^2$, and use (\ref{weakheatb}) and
Ito's lemma to derive a semimartingale decomposition for
$\int_0^t\int|\tilde u(s,x)|\Psi_s(x)dxds$, where $\Psi_s(x)\ge 0$ is a smooth
test function.  This involves approximating $|\tilde u(s,x)|$ by
$\psi_n(\langle\tilde u_s,\Phi_m(\cdot-x)\rangle)$ as $m,n\to\infty$, where
$\{\psi_n\}$ are smooth functions approximating the absolute value function
as in \cite{YW71}, and $\{\Phi_m\}$ is a smooth approximate identity.  In
Section~2 the martingale and standard drift terms which arise are handled in
a relatively straightforward manner in a general setting including that of both
Theorems~\ref{thm:uniquekbounded} and \ref{thm:uniquekunbounded} (see
Lemma~\ref{lem:conv}).  Here we may let $m,n\to\infty$ in any manner.  The
problematic term, called $I_3^{m,n}$ below, is the one arising from the
$\psi_n''/2$ term in using Ito's lemma and so will involve the quadratic
variation of the martingale term.  In the context of the Yamada-Watanabe proof,
it is the one which leads to the local time at
$0$ of the difference of two solutions to the sde, $L^0_t(X^1-X^2)$.  There,
this term is shown to be $0$ using the modulus of continuity of $\sigma$ and the
regularity of the sample paths of the solutions (the latter implicitly as one
needs the stochastic calculus associated with continuous semimartingales). 

In Section~3, $I_3^{m,n}$ is shown to approach $0$ if we first let
$m\to\infty$ and then $n\to\infty$ in the simpler context of
Theorem~\ref{thm:uniquekbounded}.  This leads to 
\begin{equation}
\label{heatinequ}
\int\IE(|\tilde u(t,x)|)\Psi_t(x)dx\le \int_0^t\int \IE(|\tilde
u(s,x)|){1\over 2}\Delta \Psi_s(x)+\dot\Psi_s(x)|dxds,
\end{equation}
from which $\tilde u=0$ follows easily by taking $\Psi_s(x)=\int
p_{t-s}(y-x)\phi(x)dx$. We feel the ease of this argument is partly related to
the greater path regularity $\tilde u$ in this context--it is H\"older
continuous in space with index $1-\epsilon$ and in time with index ${1\over
2}-\epsilon$ by results of Sanz-Sol\'e and Sarra (see \cite{SS02} and
Lemma~\ref{thm:SanzSoleSarra} below).  

In Section 4 we complete the proof of
Theorem~\ref{thm:uniquekunbounded} by showing $\lim_{n\to\infty}
I_3^{m_n,n}=0$ for a judicious choice of $m_n$, which again leads to
(\ref{heatinequ}). In this setting $\tilde u(t,x)$ is only H\"older continuous
of index ${1-\alpha/2\over 2}-\epsilon$ in time and $1-{\alpha\over 2}-\epsilon$
in space (see Lemma~\ref{thm:SanzSoleSarra} or \cite{SS02}) and this additional
irregularity makes the argument more involved.  In the Yamada-Watanabe context,
the key fact that $L^0_t(X^1-X^2)=0$ reflects the fact that the solutions must
separate ``slowly" if they do so at all.  In our setting we will argue along
similar lines by showing that $\tilde u(t,x)$ is more regular in $(t,x)$ at
small values of $\tilde u(t,x)$, i.e., when the solutions are close (see
Theorem~\ref{thm:smallHoelder}). For example, they will be H\"older of index
${1-\alpha/2\over 1-\gamma}\wedge 1-\epsilon$ in space near
space-time points where $\tilde u$ is sufficiently small
(see Corollary~\ref{cor:ubnd}).  Theorem~\ref{thm:smallHoelder} is proved in
Section 5 and is the key to the proof of Theorem~\ref{thm:uniquekunbounded}
which is completed in Section~4.  This improved modulus of continuity result may
be of independent interest. In fact a similar result to
Theorem~\ref{thm:smallHoelder} was derived independently by Mueller and Tribe
in the context of white noise, in their ongoing work on the zero set of
solutions to (\ref{sbm}).  The continuity results of Sanz-Sol\'e and  Sarra
(\cite{SS02}) and the factorization method they use (see \cite{DKZ87}), play a critical role
in the proof of Theorem~\ref{thm:smallHoelder} in our colored noise setting. 
Section~6 is an Appendix including the proofs of the weak existence theorem
(Theorem~\ref{thm:nonLipschitzexistence}) and the required moment estimates
(Proposition \ref{prop}).   

\medskip
\paragraph{Acknowledgements.}
One of us [EP] thanks Yongjin Wang for enjoyable
discussions which helped in the proof of an earlier version of
Theorem \ref{thm:smallHoelder}. The first author would like to 
express his gratitude for the
opportunity to visit  Weierstrass Institute for
Applied Analysis and Stochastics (Germany) and 
University of British Columbia (Canada) where this research was partially done.
The third author would like to express her gratitude for the opportunity
to research at the University of British Columbia (Canada) and for the
opportunity to visit the Technion (Israel).
                                                      
\section{Some auxiliary results}
\label{section:auxiliary}

Let $\rho$ be as in (\ref{rhocond}). An elementary argument shows that
$\int_{0+}(\rho(x)+\sqrt x)^{-2}dx=+\infty$ (e.g. consider
$\liminf_{x\downarrow 0}\rho^{-2}(x)x\ge 1$ and  $\liminf_{x\downarrow
0}\rho^{-2}(x)x< 1$ separately).  As we will be using $\rho$ as a
modulus of continuity (see (\ref{cond:Yamada})) we may replace $\rho$
with
$\rho(x)+\sqrt x$ and so assume
\begin{equation}\label{rhobnd}
\rho(x)\ge \sqrt x.
\end{equation}
As in the proof of Yamada and Watanabe
\cite{YW71}, we may define  a sequence of functions $\phi_{n}$ in the
following way.  First, let $a_{n} \downarrow 0$ be a strictly
decreasing sequence such that $a_{0}=1$, and
\begin{equation}
\label{acond}
\int_{a_{n}}^{a_{n-1}} \rho^{-2}(x) dx = n.
\end{equation}
Second, we define functions $\psi_{n} \in C^{\infty}_{c}(\IR)$ such that
$supp(\psi_{n}) \subset  (a_{n}, a_{n-1})$, and that
\begin{equation}
\label{psicond}
0 \leq \psi_{n}(x) \leq \frac{2 \rho^{-2}(x)}{n}\leq {2\over nx}
\quad \mbox{ for all $x \in \IR$ as well as } \quad \int_{a_{n}}^{a_{n-1}} \psi_{n}(x) dx =1.
\end{equation}
Finally, set
\begin{equation}
\label{def:phi}
\phi_{n}(x) = \int_{0}^{|x|} \int_{0}^{y} \psi_{n}(z) dz dy.
\end{equation}
From this it is easy to see that $\phi_{n}(x) \uparrow |x|$ uniformly in $x\geq 0.$
Note that each $\psi_{n}$ and thus also each $\phi_{n}$ is
identically zero in a neighborhood of zero. This implies that
$\phi_{n} \in C^{\infty}(\IR)$ despite the absolute value in its
definition. We have
\begin{eqnarray}
\label{phidiff1}
\phi_{n}'(x) &=& \sgn(x) \int_{0}^{|x|}  \psi_{n}(y) dy,\\
\label{phidiff2}
\phi_{n}''(x) &=&  \psi_{n}(|x|).
\end{eqnarray}
Thus, $|\phi_{n}'(x)| \leq 1$, and
$\int \phi_{n}''(x) h(x) dx \rightarrow h(0)$ for any function $h$ which
is continuous at zero.

\smallskip
Now let $u^{1}$ and $u^{2}$ be two solutions of (\ref{heatecolint}) with sample
paths in $C(\IR_+,C_{tem})$ a.s.,
 with the same
initial condition, $u^{1}(0)=u^{2}(0)=u_0\in C_{tem}$, and the same noise $W$
in either the setting of Theorem \ref{thm:uniquekbounded} or Theorem
\ref{thm:uniquekunbounded}.  We proceed assuming Proposition~\ref{prop} which
will be derived in the Appendix.  Define $\tilde{u} \equiv u^{1} - u^{2}.$ Let
$\Phi
\in C_{c}^{\infty}(\IR^{d})$ be a positive function with
$supp(\Phi) \subset B(0,1)$ (the open ball centered at 0 with radius 1)
such that $\int_{\IR^{d}}\Phi(x) dx=1$ and set
$\Phi_{x}^{m}( y) = m^{d}\Phi(m(x - y))$.
\noindent
Let $\la \cdot, \cdot \ra$ denote the scalar product on $L^{2}(\IR^{d}).$
By applying It\^{o}'s Formula to the semimartingale
$\la \tilde{u}_{t}, \Phi^{m}_{x} \ra$ of (\ref{weakheatb})
it follows that
\begin{eqnarray*}
& &\phi_{n}(\la \tilde{u}_{t}, \Phi_{x}^{m} \ra)\\
&=& \int_{0}^{t} \int_{\IR^{d}} \phi_{n}'(\la \tilde{u}_{s}, \Phi_{x}^{m} \ra)
\left( \sigma(u^{1}(s,y)) - \sigma(u^{2}(s,y)) \right)
\Phi_{x}^{m}(y)  W(dy ds)\\
& &+\int_{0}^{t} \phi_{n}'(\la \tilde{u}_{s}, \Phi_{x}^{m} \ra)
\la \tilde{u}_{s}, \frac{1}{2} \Delta \Phi_{x}^{m} \ra ds\\
& &+ \frac{1}{2}
\int_{0}^{t} \int_{\IR^{2d}}
\psi_{n}(|\la \tilde{u}_{s}, \Phi_{x}^{m} \ra|)
\left( \sigma(u^{1}(s,y)) - \sigma(u^{2}(s,y)) \right)
\left( \sigma(u^{1}(s,z)) - \sigma(u^{2}(s,z)) \right)\\
& &\phantom{AAAAAAAAAAAAAAAAAAAAAAAAAAAAAAAAAAA}
\Phi_{x}^{m}(y)\Phi_{x}^{m}(z) k(y,z) dy dz ds.
\end{eqnarray*}
We integrate this function of $x$ against another non-negative test
function
$\Psi \in C^{\infty}_{c}([0,t]\times\IR^{d})$.  Assume
$\Gamma\equiv\{x:\Psi_s(x)>0\ \exists s\le t\}
\subset B(0,K)$  for  some $K>0.$ We then obtain by the
classical and stochastic version of Fubini's Theorem, and arguing as
in the proof of Proposition II.5.7 of \cite{P02} to handle the time
dependence in $\psi$, that for any $t\geq 0t$,
\begin{eqnarray}
\label{eq:Iparts}
& &\left\la \phi_{n}(\la \tilde{u}_{t}, \Phi_{.}^{m} \ra), \Psi_t
\right\ra\\
\nonumber
&=& \int_{0}^{t}\int_{\IR^{d}}
\la \phi_{n}'(\la \tilde{u}_{s}, \Phi_{\cdot}^{m} \ra) \Phi_{\cdot}^{m}(y)
 , \Psi_s \ra \left( \sigma(u^{1}(s,y)) - \sigma(u^{2}(s,y)) \right)
W(dy,ds)\\
\nonumber
& &+ \int_{0}^{t} \la \phi_{n}'(\la \tilde{u}_{s}, \Phi_{.}^{m} \ra)
\la \tilde{u}_{s}, \frac{1}{2} \Delta \Phi_{.}^{m} \ra, \Psi_s \ra ds\\
\nonumber
& &+ \frac{1}{2}
\int_{0}^{t} \int_{\IR^{3d}}
\psi_{n}(|\la \tilde{u}_{s}, \Phi_{x}^{m} \ra|)
\left( \sigma(u^{1}(s,y)) - \sigma(u^{2}(s,y)) \right)
\left( \sigma(u^{1}(s,z)) - \sigma(u^{2}(s,z)) \right)\\
\nonumber
& &\phantom{AAAAAAAA}
\Phi_{x}^{m}(y)\Phi_{x}^{m}(z) k(y,z)
dy dz \Psi_s(x) dx ds+\int_0^t\la\phi_n(\la\tilde
u_s,\Phi^m_\cdot\ra) ,\dot\Psi_s\ra\,ds\\
\nonumber
&\equiv& I_{1}^{m,n}(t) + I_{2}^{m,n}(t) + I_{3}^{m,n}(t)+I_{4}^{m,n}(t).
\end{eqnarray}

\noindent

We need a calculus lemma.  
For $f\in C^2(\IR^d)$, let $\Vert
D^2f\Vert_\infty=\max_{i}\Vert{\partial^2 f\over \partial
x_i^2}\Vert_\infty$.

\begin{lemma}\label{calc} Let $f\in C_c^2(\IR^d)$ be non-negative and not
identically zero.  Then
$$\sup\Bigl\{\Bigl({{\partial f\over \partial x_i}(x)}\Bigr)^2
f(x)^{-1}:f(x)>0\Bigr\}\le 2\Vert D^2f\Vert_\infty.$$
\end{lemma}

\begin{proof} Assume first $d=1$.  Choose $x$ so that
$f(x)|f'(x)|>0$.  Without loss of generality assume $f'(x)>0$.  Let
$$x_1=\sup\{x'<x:f'(x')=0\}\in(-\infty, x).$$
By the Cauchy (or generalized mean value) theorem there is an
$x_2\in(x_1,x)$ so that
$$(f'(x)^2-f'(x_1)^2)f'(x_2)=(f(x)-f(x_1)){d((f')^2)\over dx}(x_2)
$$
and, as $f'(x_2)>0$, we get
$$f'(x)^2=(f(x)-f(x_1))2f''(x_2).
$$
Since $f$ is strictly increasing on $(x_1,x)$, and $f(x_1)\ge 0$,
$${f'(x)^2\over f(x)}\le {f'(x)^2\over f(x)-f(x_1)}\le 2\Vert
f''\Vert_\infty.
$$

For the $d$-dimensional case, assume $x$ satisfies $f(x)>0$ and let
$e_i$ be the $i$th unit basis vector. Now apply the one-dimensional
result to
$g(t)=f(x+te_i)$,
$t\in\IR$, at $t=0$.
\end{proof}

We now consider the expectation of
expression (\ref{eq:Iparts}) stopped at a stopping time $T,$
that we will specify later on.
For all the terms except $I_3^{m,n}$ we can give a unified treatment for the
settings of both Theorems~\ref{thm:uniquekunbounded} and
\ref{thm:uniquekbounded}.
\begin{lemma}\label{lem:conv} For any stopping time $T$ and constant $t\ge0$ we
have:\hfil\break
\noindent(a) \begin{equation}\label{eq:I_1}
\IE(I_1^{m,n}(t\wedge T))=0\hbox{ for all }m,n.
\end{equation}
\noindent(b) \begin{equation}
\label{eq:limI_2}
\limsup_{m,n\rightarrow \infty} \IE( I_{2}^{m,n}(t \wedge T))
\le  \IE \Big( \int_{0}^{t \wedge T} \int_{\IR}
|\tilde{u}(s,x)| {1\over 2} \Delta\Psi_s(x)dx ds \Big).
\end{equation}

\noindent(c)
\begin{equation}
\label{eq:limI_4}
\lim_{m,n\rightarrow \infty} \IE(I_4^{m,n}(t\wedge
T))=\IE\Bigl(\int_0^{t\wedge T}|\tilde u(s,x)|\dot\Psi_s(x)\,ds\Bigr).
\end{equation}

\end{lemma}

\begin{proof}
(a) Let $g_{m,n}(s,y)=\langle\phi_n'(\langle \tilde
u_s,\Phi^m_\cdot\rangle)\Phi^m_\cdot(y),\Psi_s\rangle$.  Note first that
$I_1^{m,n}(t\wedge T)$ is a continuous local martingale with square function
\begin{eqnarray}
\nonumber
\langle I_1^{m,n}\rangle_{t\wedge T}&=&\int_0^{t\wedge
T}\int\int g_{m,n}(s,y)g_{m,n}(s,z)(\sigma(u^1(s,y))-\sigma(u^2(s,y)))\\
\nonumber& &\phantom{\int_0^{t\wedge
T}\int\int}\times(\sigma(u^1(s,z))-\sigma(u^2(s,z)))k(y,z)dydzds\\
\nonumber &\le& C\int_0^{t\wedge T}\int\int
|g_{m,n}(s,y)||g_{m,n}(s,z)|(|u^1(s,y)|+|u^2(s,y)|+1)\\
\nonumber& &\phantom{C\int_0^{t\wedge T}\int\int
|g}\times(|u^1(s,z)|+|u^2(s,z)|+1)
(|z-y|^{-\alpha}+1)dydzds.
\end{eqnarray}
An easy calculation shows that $|g_{m,n}(s,y)|\le \Vert \Psi\Vert_\infty
1(|y|\le K+1)$.  Now use H\"older's inequality and (\ref{eq:uniformmoments2a})
to conclude that 
\[\IE(\langle I_1^{m,n}\rangle_{t\wedge T})\le C\int_0^t\int\int 1(|y|\le
K+1)1(|z|\le K+1)(|y-z|^{-\alpha}+1)dydzds<\infty\ \forall t>0.\]
This shows $I_1^{m,n}(t\wedge T)$ is a square integrable martingale and so has
mean $0$, as required.  

(b) In order to rewrite $I_{2}^{m,n}$
we note that both $\phi_{n}'(\la \tilde{u}_{s}, \Phi_{.}^{m} \ra)$
as well as $\la \tilde{u}_{s}, \frac{1}{2} \Delta \Phi_{.}^{m} \ra$
are in $C^{\infty}(\IR^d)$ a.s.
This follows from the infinite differentiability of the test functions
$\phi_{n}$ and $\Phi$ and from (\ref{eq:uniformmoments2a}).
Denote by $\Delta_{x}$ the Laplacian acting with respect to $x.$
Since $\tilde{u}_{s}$ is locally integrable
and $\Phi$ smooth we have for $|x|\le K$,
\begin{equation}
\label{Deltaex}
\int_{\IR^{d}} \tilde{u}(s,y)  \frac{1}{2} \Delta_{y}\Phi^{m}(x-y) dy
= \int_{\IR^{d}} \tilde{u}(s,y)  \frac{1}{2} \Delta_{x}\Phi^{m}(x-y) dy
= \frac{1}{2} \Delta_{x} \int_{\IR^{d}} \tilde{u}(s,y) \Phi^{m}(x-y) dy,
\end{equation}
for all $m.$
This implies for any $t\geq 0$,
\begin{eqnarray}
\nonumber
I_{2}^{m,n}(t)&=&
\int_{0}^{t}\int_{\IR^{d}}\phi_{n}'(\la \tilde{u}_{s}, \Phi_{x}^{m} \ra)
  \frac{1}{2} \Delta_{x}\left(\la \tilde{u}_{s},
\Phi_{x}^{m} \ra\right) \Psi_s(x) dx ds\\
\nonumber &=&
-\sum_{i=1}^{d} \frac{1}{2}
\int_{0}^{t}\int_{\IR^{d}} \frac{\partial}{\partial x_{i}}
\left(\phi_{n}'(\la \tilde{u}_{s}, \Phi_{x}^{m} \ra)\right)
\frac{\partial}{\partial x_{i}}
\left( \la \tilde{u}_{s}, \Phi_{x}^{m} \ra \right) \Psi_s(x) dx ds\\
\nonumber & &- \sum_{i=1}^{d} \frac{1}{2}
\int_{0}^{t}\int_{\IR^{d}} \phi_{n}'(\la \tilde{u}_{s}, \Phi_{x}^{m} \ra)
\frac{\partial}{\partial x_{i}}
\left( \la \tilde{u}_{s}, \Phi_{x}^{m} \ra \right)
\frac{\partial}{\partial x_{i}}\Psi_s(x) dx ds\\
\nonumber &=&-   \sum_{i=1}^{d}  \frac{1}{2}\int_{0}^{t}\int_{\IR^{d}}
\psi_{n}(|\la \tilde{u}_{s}, \Phi_{x}^{m} \ra|)
\left( \frac{\partial}{\partial x_{i}}
\la \tilde{u}_{s}, \Phi_{x}^{m} \ra \right)^{2} \Psi_s(x) dx ds\\
\nonumber & & - \sum_{i=1}^{d}  \frac{1}{2}
\int_{0}^{t}\int_{\IR^{d}}\phi_{n}'(\la \tilde{u}_{s}, \Phi_{x}^{m} \ra)
\frac{\partial}{\partial x_{i}}
\left( \la \tilde{u}_{s}, \Phi_{x}^{m} \ra  \right)
\frac{\partial}{\partial x_{i}}\Psi_s(x) dx ds
\end{eqnarray}
\begin{eqnarray}
\nonumber
&=& -   \sum_{i=1}^{d}  \frac{1}{2}\int_{0}^{t}\int_{\IR^{d}}
\psi_{n}(|\la \tilde{u}_{s}, \Phi_{x}^{m} \ra|)
\left( \frac{\partial}{\partial x_{i}}
\la \tilde{u}_{s}, \Phi_{x}^{m} \ra \right)^{2} \Psi_s(x) dx ds\\
\nonumber
& &+\sum_{i=1}^{d}  \frac{1}{2}\int_{0}^{t}\int_{\IR^{d}}
\psi_{n}(\la \tilde{u}_{s},\Phi_{x}^{m}\ra)\frac{\partial}{\partial x_{i}}
(\la \tilde{u}_{s},\Phi_{x}^{m}\ra)
\la \tilde{u}_{s},\Phi_{x}^{m}\ra \frac{\partial}{\partial x_{i}} \Psi_s(x)
dx ds\\
\nonumber
& &+ \int_{0}^{t}\int_{\IR^{d}}
\phi_{n}'(\la \tilde{u}_{s},\Phi_{x}^{m}\ra)
\la \tilde{u}_{s},\Phi_{x}^{m}\ra
\frac{1}{2} \Delta \Psi_s(x)dx ds\\
\label{secondterm}
&=& \int_{0}^{t} I_{2,1}^{m,n}(s) + I_{2,2}^{m,n}(s) + I_{2,3}^{m,n}(s) ds.
\end{eqnarray}
Above, we have used that $\phi_{n}''=\psi_{n}$ and we have repeatedly
used integration by parts, the product rule as well as the chain rule on
$\phi_{n}'(\la \tilde{u}_{s}, \Phi_{x}^{m} \ra).$
In order to deal with the various parts of $I_{2}^{m,n}$
we will first jointly consider $I_{2,1}^{m,n}$ and $I_{2,2}^{m,n}.$
For fixed $s$ and $i=1,\dots, d$ we define a.s.,
\begin{eqnarray*}
A_{i}^{s}&=& \Bigl\{ x: \left( \frac{\partial}{\partial x_i}
\la \tilde{u}_{s}, \Phi_{x}^{m} \ra \right)^{2} \Psi_s(x) \leq
\la \tilde{u}_{s}, \Phi_{x}^{m} \ra \frac{\partial}{\partial x_i}
\la \tilde{u}_{s}, \Phi_{x}^{m} \ra \frac{\partial}{\partial x_i}
\Psi_s(x) \Bigr\} \cap \{x:\Psi_s(x)>0\}\\
&=& A^{+,s}_{i}\cup A^{-,s}_{i} \cup A^{0,s}_{i},
\end{eqnarray*}
where
\begin{eqnarray*}
A^{+,s}_{i}&=& A^{s}_{i} \cap \{ \frac{\partial}{\partial x_i}
\la \tilde{u}_{s}, \Phi_{x}^{m} \ra >0\},\\
A^{-,s}_{i} &=& A^{s}_{i} \cap \{ \frac{\partial}{\partial x_i}
\la \tilde{u}_{s}, \Phi_{x}^{m} \ra <0\},\\
A^{0 ,s}_{i}&=& A^{s}_{i}  \cap \{ \frac{\partial}{\partial x_i}
\la \tilde{u}_{s}, \Phi_{x}^{m} \ra =0\}.
\end{eqnarray*}
On $A^{+,s}_{i}$ we have
\begin{equation*}
0 <\Bigl(\frac{\partial}{\partial x_i}
\la \tilde{u}_{s}, \Phi_{x}^{m} \ra \Bigr)\Psi_s(x)\leq
\la \tilde{u}_{s}, \Phi_{x}^{m} \ra \frac{\partial}{\partial x_i}
\Psi_s(x),
\end{equation*}
and therefore for any $t\geq 0$,
\begin{eqnarray*}
& &\int_0^t\int_{A^{+,s}_{i}} \psi_{n}(|\la \tilde{u}_{s}, \Phi_{x}^{m} \ra
|)
\la \tilde{u}_{s}, \Phi_{x}^{m} \ra \frac{\partial}{\partial x_i} \Psi_s(x)
\frac{\partial}{\partial x_i}
\la \tilde{u}_{s}, \Phi_{x}^{m} \ra dx\,ds\\
&\leq& \int_0^t\int_{A^{+,s}_{i}} \psi_{n}(|\la \tilde{u}_{s}, \Phi_{x}^{m}
\ra |)
\la \tilde{u}_{s}, \Phi_{x}^{m} \ra^{2} \frac{(\frac{\partial}
{\partial x_i} \Psi_s(x))^{2}}{\Psi_s(x)} dx\,ds\\
&\leq& \int_0^t\int_{A^{+,s}_{i}} \frac{2}{n}
1_{\{a_{n-1} \leq |\la \tilde{u}_{s}, \Phi_{x}^{m} \ra| \leq a_{n} \}}
|\la \tilde{u}_{s}, \Phi_{x}^{m} \ra| \frac{(\frac{\partial}
{\partial x_i} \Psi_s(x))^{2}}{\Psi_s(x)} dx\,ds\quad\hbox{by (\ref{psicond})}
\\ &\leq& \frac{2a_{n}}{n} \int_0^t\int_{\IR^d} 1(\Psi_s(x)>0)
\frac{(\frac{\partial}{\partial x_i} \Psi(x))^{2}}{\Psi_s(x)} dx\,ds\\
&\le& \frac{2a_n}{n}\int_0^t2\Vert
D^2\Psi_s\Vert_\infty\hbox{Area}(\Gamma)\,ds\equiv{2a_n\over n}C(\Psi),
\end{eqnarray*}
where Lemma~\ref{calc} is used in the last line.
Similarly, on the set $A^{-,s}_{i},$
\begin{equation*}
0 > \frac{\partial}{\partial x_i}
\la \tilde{u}_{s}, \Phi_{x}^{m} \ra \Psi_s(x) \geq
\la \tilde{u}_{s}, \Phi_{x}^{m} \ra \frac{\partial}{\partial x_i}
\Psi_s(x).
\end{equation*}
Hence, with the same calculation
\begin{eqnarray*}
& &\int_0^t\int_{A^{-,s}_{i}} \psi_{n}(|\la \tilde{u}_{s}, \Phi_{x}^{m} \ra
|)
\la \tilde{u}_{s}, \Phi_{x}^{m} \ra \frac{\partial}{\partial x_i} \Psi_s(x)
\frac{\partial}{\partial x_i}
\la \tilde{u}_{s}, \Phi_{x}^{m} \ra dx\,ds\\
&\leq& \frac{2a_{n}}{n}\int_0^t \int_{\IR^d} 1(\Psi_s(x)>0)
\frac{(\frac{\partial}{\partial x_i} \Psi_s(x))^{2}}{\Psi_s(x)} dx\,ds\\
&\le& {2a_n\over n}C(\Psi).
\end{eqnarray*}
Finally, for any $t\geq 0$,
\begin{equation*}
\int_0^t\int_{A^{0,s}_{i}} \psi_{n}(|\la \tilde{u}_{s}, \Phi_{x}^{m} \ra |)
\la \tilde{u}_{s}, \Phi_{x}^{m} \ra \frac{\partial}{\partial x_i} \Psi_s(x)
\frac{\partial}{\partial x_i}
\la \tilde{u}_{s}, \Phi_{x}^{m} \ra dx\,ds=0,
\end{equation*}
and we conclude that
\begin{equation*}
\IE( I_{2,1}^{m,n}(t \wedge T)
+ I_{2,2}^{m,n}(t \wedge T)) \leq 4C(\Psi)
 \frac{a_{n}}{n},
\end{equation*}
which tends to zero as $n \rightarrow \infty.$ For $I_{2,3}^{m,n}$
recall that $\phi_{n}'(u) u \uparrow |u|$ uniformly in $u$ as
$n \rightarrow \infty,$
and that $\la \tilde{u}_{s}, \Phi_{x}^{m} \ra$
tends to $\tilde{u}(s,x)$ as $m \rightarrow
\infty$ for all $s,x$ a.s. by the a.s. continuity of $\tilde u$.
This implies
that
$\phi_{n}'(\la \tilde{u}_{s}, \Phi_{x}^{m} \ra)
\la \tilde{u}_{s}, \Phi_{x}^{m} \ra
\rightarrow |\tilde{u}(s,x)|$
pointwise a.s.
as $m,n \rightarrow \infty,$ where it is unimportant how we take the limit.
We also have the bound
\begin{eqnarray}
\label{phinbound}
|\phi_{n}'(\la \tilde{u}_{s}, \Phi_{x}^{m} \ra)
\la \tilde{u}_{s}, \Phi_{x}^{m} \ra|
&\leq&
|\la \tilde{u}_{s}, \Phi_{x}^{m} \ra|
\leq \la |\tilde{u}_{s}|, \Phi_{x}^{m} \ra.
\end{eqnarray}
The a.s. continuity of $\tilde u$ implies a.s. convergence for all $s,x$ of
$\langle|\tilde u_s|,\Phi^m_x\rangle$ to $|\tilde u(s,x)|$ as $m\to \infty$. A
simple application of Jensen's Inequality and (\ref{eq:uniformmoments2a}) shows
that
$|\langle|\tilde u_s|,\Phi^m_x\rangle \,|$ is $L^p$
bounded on $([0,t]\times B(0,K)\times \Omega,ds\times
dx\times \IP)$ uniformly in $m.$  This implies
\begin{equation}\label{ui}
\{\la |\tilde{u}_{s}|, \Phi_{x}^{m} \ra
:m\}\hbox{ is uniformly integrable
on }([0,t]\times B(0,K)\times \Omega).
\end{equation}
 and so gives uniform
integrability of $
\{|\phi_{n}'(\la \tilde{u}_{s}, \Phi_{x}^{m} \ra)
\la \tilde{u}_{s}, \Phi_{x}^{m} \ra|:m,n\}$ by our earlier bound
(\ref{phinbound}).
This implies
 $$\lim_{m,n\to\infty}\IE( I_{2,3}^{m,n}(t\wedge
T))=\IE\Bigl(\int_0^{t\wedge T}\int|\tilde u(s,x)|{1\over
2}\Delta \Psi_s(x)dx\,ds\Bigr).$$  Collecting the pieces, we have shown that
(\ref{eq:limI_2}) holds.

\medskip

\noindent(c) As in the above argument we have
\begin{equation}\label{absconv}
\phi_n(\la\tilde u_s,\Phi_x^m\ra)\to|\tilde
u(s,x)|\hbox{ as }m,n\to\infty\hbox{ a.s. for all }x\hbox{ and all }s\le t.
\end{equation}
The uniform integrability in (\ref{ui})
and the bound $\phi_n(\la\tilde u_s,\Phi^m_x\ra)\le\la|\tilde
u_s|,\Phi_x^m\ra$ imply
$$\{\phi_n(\la\tilde u_s,\Phi_x^m\ra:n,m\}\hbox{ is uniformly integrable
on }[0,t]\times B(0,K)\times\Omega.$$
Therefore the result now follows from the above convergence and the
bound $$|\dot\Psi_s(x)|\le
C1(|x|\le K).$$ \end{proof}

\section{Proof of Theorem \ref{thm:uniquekbounded}}
\label{section:proofkbounded}
\noindent
Here, we let $T=t$ be deterministic.  Given the results from 
Section~\ref{section:auxiliary} it now remains to estimate
 $\IE( I_{3}^{m,n}(t)).$ We will then let
$m \rightarrow \infty$ before letting $n \rightarrow \infty.$
By the boundedness of the correlation
kernel $k$ and Jensen's Inequality, $I_{3}^{m,n}(t)$ is bounded by
\begin{eqnarray*}
& &\frac{1}{2} ||k||_{\infty} \int_{0}^{t} \int_{\IR^{d}}
 \left( \int_{\IR^{d}} \left| \sigma(u^{1}(s,y))
- \sigma(u^{2}(s,y)) \right| \Phi_{x}^{m}(y) dy \right)^{2}
\psi_{n}(|\la \tilde{u}_{s}, \Phi_{x}^{m} \ra|) \Psi_s(x) dx   ds\\
&\leq& \frac{1}{2}||k||_{\infty} \int_{0}^{t} \int_{\IR^{d}} \left(
\sigma(u^{1}(s,y))  - \sigma(u^{2}(s,y)) \right)^{2}
\left( \int_{\IR^{d}} \psi_{n}(|\la \tilde{u}_{s}, \Phi_{x}^{m} \ra|)
  \Phi_{x}^{m}(y) \Psi_s(x) dx \right) dy  ds.
\end{eqnarray*}
The integral in parentheses is bounded by a
constant, independent of $m$, is zero for all $m$ if $|y|>K+1$, and as
$m\to 0$ converges to
$ \psi_{n}(\tilde{u}(s,y)) \Psi_s(y)$  for all $(s,y)$ by
the continuity of $\tilde u$.
Our growth condition on $\sigma$ and (\ref{eq:uniformmoments2a}) imply the
integrability of 
\[\int_0^t\int(\sigma(u^1(s,y)-\sigma(u^2(s,y))^2  1_{\{|y|\le
K+1\}}\,dyds.\]
Therefore, the Dominated Convergence Theorem
implies that
\begin{eqnarray}
\nonumber
\limsup_{m \rightarrow \infty} \IE \left(  I_{3}^{m,n}(t) \right)
&\leq& \frac{1}{2}||k||_{\infty} \IE \left( \int_{0}^{t}  \la
\psi_{n}(\tilde{u}_{s})\left( \sigma(u^{1}_{s})
- \sigma(u^{2}_{s}) \right)^{2},
\Psi_s \ra ds \right)\\
\label{Phixconv}
&\leq& C(\Psi)||k||_{\infty}\frac{t}{n}  ,
\end{eqnarray}
where the last line follows by (\ref{cond:Yamada}) and (\ref{psicond}).

Return to equation (\ref{eq:Iparts}) and let first
$m \rightarrow
\infty$ and then $n \rightarrow \infty$.  Use the above and Lemma~\ref{lem:conv}
on the right-hand side, and
(\ref{absconv}) and Fatou's lemma on the left-hand side, to conclude that
\begin{equation}
\label{eq:heatGron}
\int_{\IR^{d}} \IE\left( |\tilde{u}(t,x)|\right) \Psi_t(x) dx
\leq \int_{0}^{t}\int_{\IR^{d}} \IE\left( |\tilde{u}_{s}(x)|\right)
|\frac{1}{2} \Delta \Psi_s(x)+\dot\Psi_s(x)| dx ds.
\end{equation}
Let $\{g_N\}$ be a sequence of functions in $C_c^\infty(\IR^d)$ such that
$g_N:\IR^d\to[0,1]$,
$$B(0,N)\subset\{x:g_N(x)=1\},\quad B(0,N+1)^c\subset\{x:g_N(x)=0\},$$
and
$$\sup_N[\Vert\nabla g_N\Vert_\infty+\Vert D^2g_N\Vert_\infty]\equiv
C(g)<\infty,$$
where $\nabla g_N$ denotes the gradient with respect to the spatial variables.
Now let $\phi\in C_c^\infty(\IR^d)$, and for $(s,x)\in[0,t]\times\IR^d$
set $\Psi_N(s,x)=(S_{t-s}\phi(x)) g_N(x)$.  It is then easy to check that $\Psi_N\in
C_c^\infty([0,t]\times \IR^d)$ and for $\lambda>0$ there is a
$C=C(\lambda,\phi)$ such that for all $N$
\begin{eqnarray*}
|{\Delta\over
2}\Psi_N(s,x)+\dot\Psi_N(s,x)|&=&
\left|\sum_{i=1}^{d} \frac{\partial}{\partial x_{i}}
S_{t-s}\phi(x_{i})\frac{\partial}{\partial x_{i}}g_N(x_{i})
+S_{t-s}\phi(x){\Delta\over 2}g_N(x)\right|\\
&\le&C e^{-\lambda|x|}1_{\{|x|>N\}}.
\end{eqnarray*}
Use this in (\ref{eq:heatGron}) to conclude that
$$\int_{\IR^d}\IE(|\tilde u(t,x)|)\phi(x)\,dx\le
C\int_0^t\int_{\IR^d}\IE(|\tilde u(s,x)|)e^{-\lambda
|x|}1_{\{|x|>N\}}\,dx\,ds.$$
By Proposition~\ref{prop} 
the right-hand side of the above approaches zero as $N\to\infty$ and we
see that
$$\IE\Bigl(\int_{\IR^d}|\tilde u(t,x)|dx\Bigr)=0.$$
Therefore $u^1(t)=u^2(t)$ for all $t\ge 0$ a.s. by a.s. continuity.

\section{Proof of Theorem \ref{thm:uniquekunbounded}}
\label{s:kunbounded}
We continue to use the notation of Section~2 and also assume the hypotheses of
Theorem~\ref{thm:uniquekunbounded}.  In particular $u^1$ and $u^2$ are solutions of
(\ref{heatecolint}), $\tilde u=u^1-u^2$, 
$\sigma$ is H\"older continuous with exponent
$\gamma:$
\begin{equation*}
|\sigma(u) -\sigma(v)| \leq L |u - v|^{\gamma} \hbox{ for } u,v \in \IR,
\end{equation*}
and $|{k}(x,y)|\leq c_{\ref{thm:uniquekunbounded}} [|x-y|^{-\alpha}+1]$ for some
$\alpha\in(0,1)$.  We choose
$\rho(x)=\sqrt{x}$ for our smooth approximation of the absolute
value function  throughout noting that (\ref{cond:Yamada}) is not
necessarily  satisfied for large values.
Nevertheless, we will use the test function $\phi_{n}$ and its derivatives
as defined in (\ref{def:phi}) to (\ref{phidiff2}) corresponding to this $\rho.$

Fix some $\lambda>0$ and
let $T_K=\inf\{t\ge 0: \sup_{x \in \IR^{d}} (|u^1(t,x)|+|u^2(t,x)|) e^{-\lambda |x|}>K\}
\wedge K.$ 
Note that 
\begin{eqnarray}
\label{equt:120}
T_{K}\rightarrow \infty, \;\; \;P-{\rm a.s.},\;
\end{eqnarray}
since $u^i\in C(\IR_{+}, C_{tem})$. 

Also define a metric
$d$ by
\[d((t,x),(t',x'))=\sqrt{|t-t'|}+|x-x'|, t,t'\in\IR_+, x,x'\in\IR^d,\]
and set
\begin{eqnarray*} Z_{K,N,\xi}=\{(t,x)\in\IR_+\times\IR^d:&t\le T_K,|x|\le K,
d((t,x),(\hat{t},\hat{x}))<2^{-N}\hbox{ for some } \\
&(\hat{t},\hat{x})\in[0,T_K]\times \IR^d \hbox{
satisfying }|\tilde u(\hat{t},\hat{x})|\le 2^{-N\xi}\}.
\end{eqnarray*}

We will now use the following key result on the improved H\"older continuity of
$\tilde u$ when $\tilde u$ is small.  It will be proved in Section
\ref{proof:smallHoelder}.
\begin{thm}
\label{thm:smallHoelder}
Assume the hypotheses of Theorem~\ref{thm:uniquekunbounded}, except now
allow $\gamma\in (0,1]$.  Let $u_0\in C_{tem}$ and $\tilde u=u^1-u^2$,
where $u^i$ is a solution of \eqref{heatecolint} with sample paths in
$C(\IR_+,C_{tem})$ a.s. for
$i=1,2$.
Let $\xi\in (0,1)$ satisfy
\begin{eqnarray}
\nonumber
&\exists N_\xi=N_\xi(K,\omega)\in\IN\
a.s.\hbox{ such that }\forall N\ge N_\xi,
(t,x)\in Z_{K,N,\xi}\\
\label{reghyp}
&d((t',y),(t,x))\le 2^{-N}, t,t'\le
T_K\Rightarrow |\tilde u(t,x)-\tilde u(t',y)|\le
2^{-N\xi}.
\end{eqnarray}
Let $0<\xi_1<[\xi \gamma + 1-\frac{\alpha}{2}]\wedge 1$.
Then there is an $N_{\xi_1}=N_{\xi_1}(K,\omega)\in\IN$
a.s. such that for any $N\ge N_{\xi_1}$ in $\IN$ and any $(t,x)\in Z_{K,N,\xi}$
\begin{equation}
\label{reghyp2}
d((t',y),(t,x))\le 2^{-N}, t, t'\le
T_K\Rightarrow |\tilde u(t,x)-\tilde u(t',y)|\le 2^{-N\xi_1}.
\end{equation}
Moreover there are strictly positive constants
$R,\delta, c_{\ref{Nbnd}.1},c_{\ref{Nbnd}.2}$ depending only on
$(\xi,\xi_1)$ and $N(K) \in \IN,$ which also depends on $K,$ such
that
\begin{eqnarray}
\label{Nbnd}
\IP(N_{\xi_1}\ge N)\le c_{\ref{Nbnd}.1} ( \IP(N_\xi\ge N/R)
+ K^{d+1} \exp(-c_{\ref{Nbnd}.2} 2^{N\delta}) )
\end{eqnarray}
provided that $N\geq N(K).$
\end{thm}

\noindent{\bf Remark.} Results similar to the above for white noise were
independently found by Carl Mueller and Roger Tribe in their parallel
work on level sets of solutions of SPDE's.
\medskip

\noindent Recall $\lambda>0$ is a fixed parameter used in the definition of
$T_K$. 
\begin{cor}
\label{cor:ubnd}
Assume the hypthoses of Theorem~\ref{thm:uniquekunbounded} except now allow
$\gamma\in(0,1]$.  Let $u_0$ and $\tilde u$ be as in
Theorem~\ref{thm:smallHoelder}, and $1-\frac{\alpha}{2}<\xi<
\frac{1-\frac{\alpha}{2}}{1-\gamma} \wedge 1$. There is an a.s. finite positive
random variable
$C_{\xi,K}(\omega)$ such that for any $\epsilon\in(0,1]$, $t\in[0,T_K]$ and $|x|\le
K$, if $|\tilde u(t,\hat{x})|\le \epsilon^\xi$ for some $|\hat{x}-x|\le \epsilon$, then
$|\tilde u(t,y)| \le C_{\xi,K}\epsilon^\xi$ whenever $|x-y|\le \epsilon$.  Moreover
there are strictly positive constants
$\delta,c_{\ref{Cprobbnd}.1},c_{\ref{Cprobbnd}.2}$, depending on $\xi$, and an
$r_0(K),$ which also depends on $K,$ such that
\begin{eqnarray}
\label{Cprobbnd}
\IP(C_{\xi,K}\ge r)\le c_{\ref{Cprobbnd}.1}\Bigl[ \Bigl(\frac{r-6}{(K+1)
e^{\lambda (K+1)}}\Bigr)^{-\delta} + K^{d+1}
\exp\Bigl(-c_{\ref{Cprobbnd}.2}\Bigl({r-6\over (K+1) e^{\lambda
(K+1)}}\Bigr)^\delta\Bigr)\Bigr]
\end{eqnarray}
for all $r\geq r_0(K)>6+(K+1) e^{\lambda (K+1)}.$
\end{cor}
\begin{proof}
By  Proposition~\ref{prop}(b) and the equality $\tilde u=Z^1-Z^2$, where
$Z^i(t,x)=u^i(t,x)-S_tu_0(x)$, we have (\ref{reghyp})  with
$\xi=\xi_0=\frac{1}{2}(1-\frac{\alpha}{2})$.  Indeed, $\tilde{u}$ is
uniformly H\"older continuous on compacts in $[0,\infty)\times\IR^d$ with
coefficient
$\xi$ in space and
$\frac{\xi}{2}$ in time provided that $\xi<1-\frac{\alpha}{2}.$

Inductively define $\xi_{n+1}=\Bigl[\Bigl( \xi_{n} \gamma + 1
-\frac{\alpha}{2}\Bigr)
\wedge 1\Bigr]\Bigl(1-{1\over n+3}\Bigr)$
so that
$\xi_n\uparrow \frac{1-\frac{\alpha}{2}}{1-\gamma}\wedge 1$.
Fix $n_0$ so that $\xi_{n_0}\ge \xi>\xi_{n_0-1}$.  Apply
Theorem~\ref{thm:smallHoelder} inductively $n_0$ times to get (\ref{reghyp}) for
$\xi_{n_0-1}$ and hence (\ref{reghyp2}) with $\xi_1=\xi_{n_0}$.  

First consider $\epsilon\le 2^{-N_{\xi_{n_0}}}$.  Choose $N\in\IN$ so that
$2^{-N-1}<\epsilon\le 2^{-N}$ $(N\ge N_{\xi_{n_0}})$, and assume $t\le T_K$, $|x|\le
K$ and
$|\tilde u(t,\hat{x})|\le \epsilon^\xi\le 2^{-N\xi}\le 2^{-N\xi_{n_0-1}}$ for some
$|\hat{x}-x|\le \epsilon\le 2^{-N}$. Then $(t,x)\in Z_{K,N,\xi_{n_0-1}}$.  Therefore
(\ref{reghyp2}) with $\xi_1=\xi_{n_0}$ implies that if $|y-x|\le \epsilon\le
2^{-N}$, then
\begin{eqnarray*}
|\tilde u(t,y)|&\le&|\tilde u(t,\hat{x})|+|\tilde u(t,\hat{x})-\tilde u(t,x)|
+|\tilde u(t,x)-\tilde u(t,y)|\\
&\le &2^{-N\xi}+2\cdot 2^{-N\xi_{n_0}}
\le 3\cdot 2^{-N\xi}\le 3(2\epsilon)^\xi\le 6\epsilon^\xi.
\end{eqnarray*}
For $\epsilon>2^{-N_{\xi_{n_0}}}$, we have for $(t,x)$
and $(t,y)$ as in the corollary,
\[|\tilde u(t,y)|\le (K+1) e^{\lambda (K+1)}
\le (K+1) e^{\lambda (K+1)} 2^{N_{\xi_{n_0}}}\epsilon^\xi. \]
This gives the conclusion with $C_{\xi,K}=(K+1) e^{\lambda (K+1)} 2^{N_{\xi_{n_0}}}+6$.
A short calculation and (\ref{Nbnd}) now imply that there are strictly positive
constants
$\tilde{R},\tilde{\delta},{c}_{\ref{Cprobbnd1}.1},{c}_{\ref{Cprobbnd1}.2}$,
depending on $\xi$ and $K$, such that
\begin{eqnarray}
\nonumber
\IP(C_{\xi,K}\ge r)&\le&
{c}_{\ref{Cprobbnd1}.1}\Bigl[\IP\Bigl(N_{\frac{1}{2}(1-\frac{\alpha}{2}) }\ge
\frac{1}{\tilde{R}} \log_2 \Bigl({r-6\over (K+1) e^{\lambda (K+1)}}\Bigr)
\Bigr)\\
&&\phantom{\tilde{c}_{\ref{Cprobbnd1}.1}\Bigl[\IP\Bigl(N_{\frac{1}{2}(1-\frac{\alpha}{2})
}\ge \frac{1}{R} \log_2}
\label{Cprobbnd1}
+ K^{d+1} \exp\Bigl(-{c}_{\ref{Cprobbnd1}.2}\Bigl({r-6\over (K+1)
e^{\lambda (K+1)}}\Bigr)^{\tilde{\delta}}\Bigr)\Bigr]
\end{eqnarray}
for all $r\geq r_0(K).
$
The usual Kolmogorov continuity proof applied to 
(\ref{eq:Kolmogorov}) with
$\tilde{u}=Z^1-Z^2$ in place of $Z$ (and $\xi=\frac{1}{2}(1-\frac{\alpha}{2})$)
shows there are $\tilde{\epsilon},\tilde{c}_{3} >0$
such that
$$\IP(N_{\frac{1}{2}(1-\frac{\alpha}{2})}\geq M)
\leq \tilde{c}_{3} 2^{-M \tilde{\epsilon}}$$
for all $M \in \IR.$
Thus, (\ref{Cprobbnd}) follows from (\ref{Cprobbnd1}).
\end{proof}

\medskip

Now fix $\alpha, \gamma$ satisfying the conditions of
Theorem \ref{thm:uniquekunbounded}, so $\alpha<(2\gamma-1)$
and notice that since $1\ge \gamma > \frac{1}{2}$ this implies that
$\frac{1-\frac{\alpha}{2}}{1-\gamma} > 1.$ Hence, we can choose $\xi
\in (0,1)$ such that
\begin{eqnarray}\label{alphbnd}
\alpha<\xi(2\gamma-1)
\end{eqnarray}
and $1-\frac{\alpha}{2} < \xi < \frac{1-\frac{\alpha}{2} }{1-\gamma}
\wedge 1.$ This means that $\xi$ satisfies the conditions of Corollary
\ref{cor:ubnd}.

We return to the setting and notation in Section 2. In particular $\Psi\in
C_c^\infty([0,t]\times\IR^d)$ with $\Gamma=\{x:\Psi_s(x)>0\
\exists s\le t\}\subset B(0,K)$.  Recall Lemma~\ref{lem:conv} is
valid in the setting of Theorem~\ref{thm:uniquekunbounded}.

Let $m^{(n)}:=a_{n-1}^{-\frac{1}{\xi}}$. Note that $m^{(n)}\ge
1$ for all $n$.
We set $c_{0}(K):= r_0(K) \vee K^{2}
e^{\lambda K}$ (where $r_0(K)$ is chosen as in Corollary \ref{cor:ubnd})
and define the stopping time
\begin{eqnarray*}
& &T_{\xi,K}= \inf\{ t \geq 0: t>T_{K} \text { or }t \leq T_{K} \text{ and there exist }
\epsilon \in (0,1], \hat{x},x,y \in \IR \text{ with } \\
& & \phantom{AAAA}|x|\le K, |\tilde{u}(t,\hat{x})| \leq \epsilon^{\xi},
|x-\hat{x}| \leq \epsilon, |x-y|\leq \epsilon \text{ such that } |\tilde{u}(t,y)|> c_{0}(K)
\epsilon^{\xi} \}.
\end{eqnarray*}
Assuming our filtration is completed as usual, $T_{\xi,K}$ is a stopping time 
by the standard projection argument.
Note that for any $t\geq 0,$ by Corollary \ref{cor:ubnd},
\begin{eqnarray}
\nonumber
\IP( T_{\xi,K} \leq t) &\leq& \IP( T_K \leq t) + \IP(C_{\xi,K} >c_0(K)) \\
\label{eq:stoptbnd}
&\leq& \IP( T_K \leq t)
+ c_{\ref{Cprobbnd}.1}\Bigl[ \Bigl(\frac{K^{2}e^{\lambda K}-6}{(K+1)
e^{\lambda (K+1)}}\Bigr)^{-\delta}\\
& &\phantom{\IP( T_{\xi,K} \leq t) \leq  }
+K^{d+1} \exp\Bigl(-c_{\ref{Cprobbnd}.2}\Bigl({K^{2}e^{\lambda K}-6\over (K+1)
e^{\lambda (K+1)}}\Bigr)^\delta\Bigr)\Bigr]
\end{eqnarray}
which tends to zero as $K \rightarrow \infty$ 
due to~(\ref{equt:120}). 

With this set-up we can show the following lemma:
\begin{lemma}
\label{lem:ubnd}
For all $x \in \Gamma$ and $s\in[0,T_{\xi,K}]$, if
$ |\la \tilde{u}_{s}, \Phi_{x}^{m^{(n)}} \ra| \leq a_{n-1}$ then
\begin{equation*}
\sup_{y \in B(x, \frac{1}{m^{(n)}})} |\tilde{u}(s,y)| \leq c_{0}(K) a_{n-1}.
\end{equation*}
\end{lemma}
\begin{proof}
Since $ |\la \tilde{u}_{s}, \Phi_{x}^{m^{(n)}} \ra| \leq a_{n-1}$ and $\tilde
u_s(\cdot)$ is continuous there exists a $\hat{x} \in B(x, \frac{1}{m^{(n)}})$ such
that
$|\tilde{u}(s,\hat{x})| \leq a_{n-1}.$
Apply the definition of the stopping time
with $\epsilon=1/m^{(n)}\in (0,1]$ and so $\epsilon^\xi=a_{n-1}$ to obtain the required
bound.
\end{proof}

\smallskip

Next, we bound
$|I_{3}^{m^{(n)},n}|$ using the H\"older continuity of
$\sigma,$ as well  as the definition of $\psi_{n}$.  If
$|\sigma(x)-\sigma(y)|\le L|x-y|^\gamma$, then
\begin{eqnarray*}
|I_{3}^{m^{(n)},n}(t \wedge T_{\xi,K})|
&\leq& {c_{\ref{def:singulartildek}} L^2\over n}
\int_{0}^{t\wedge T_{\xi,K}} \int_{\IR^{3}}
1_{\{a_{n} \leq |\la \tilde{u}_{s}, \Phi_{x}^{m^{(n)}} \ra| \leq a_{n-1}\} } \,
a_{n}^{-1} |\tilde{u}_{s}(y)|^{\gamma} \, |\tilde{u}_{s}(z)|^{\gamma}   \\
& &\phantom{AAAAAAAAA}\cdot
\Phi_{x}^{m^{(n)}}(y)\Phi_{x}^{m^{(n)}}(z)[|y-z|^{-\alpha}+1] dy dz \Psi_s(x) dx
ds.
\end{eqnarray*}
Now set $\Gamma^{1} = \{x \in \IR^{d}, d(x,\Gamma) < 1 \}.$ Since
$\Phi(x) \leq C 1_{B(0,1)}(x)$ and
\[1_{B(0,1)}(m^{(n)}(x-y) ) \cdot 1_{B(0,1)}( m^{(n)}(x-z) )
\leq 1_{B(0,1)}(m^{(n)} (x-y) ) \cdot 1_{B(0,1)}({1\over 2}m^{(n)}(y-z) ),\]
we obtain from Lemma \ref{lem:ubnd}
\begin{eqnarray*}
& &|I^{m^{(n)},n}_{3}(t \wedge T_{\xi,K})| \\
&\leq& c_{\ref{def:singulartildek}} L^2c_{0}(K)^{2\gamma} \frac{a_{n-1}^{2
\gamma}}{na_{n}}
\int_{0}^{t\wedge T_{\xi,K}} \int_{\IR^{3d}}
1_{\{a_{n} \leq |\la \tilde{u}_{s}, \Phi_{x}^{m^{(n)}} \ra| \leq a_{n-1}\} } \,\\
& &\phantom{AAAAAAAAAAAAAAAAAAAA}
\cdot \Phi_{x}^{m^{(n)}}(y)\Phi_{x}^{m^{(n)}}(z) [|y-z|^{-\alpha}+1]
dy dz \Psi_s(x) dx ds
\end{eqnarray*}
\begin{eqnarray*}
&\leq& \frac{c_{\ref{def:singulartildek}} L^2 
||\Psi||_{\infty} c_{0}(K)^{2\gamma}}{n} \frac{a_{n-1}^{2 \gamma}}{a_{n}}
\int_{0}^{t\wedge T_{\xi,K}} \int_{\Gamma^{1} \times \Gamma^{1}}
\left( \int_{\Gamma}
\Phi_{x}^{m^{(n)}}(y)\Phi_{x}^{m^{(n)}}(z) dx \right)[|y-z|^{-\alpha}+1]
dy dz  ds\\
&\leq &\frac{c_{\ref{def:singulartildek}} L^2 
||\Psi||_{\infty} c_{0}(K)^{2\gamma}t}{n} \frac{a_{n-1}^{2 \gamma}}{a_{n}}
 \int_{\Gamma^{1} \times \Gamma^{1}}
(m^{(n)})^{d} 1_{B(0,1)}({1\over 2}m^{(n)}(y-z) )[|y-z|^{-\alpha}+1]
dy dz\\
&\leq& \frac{C(c_{\ref{def:singulartildek}}, L, \Psi, \Phi) c_{0}(K)^{2\gamma} t}{n} 
\frac{a_{n-1}^{2
\gamma}}{a_{n}} [(m^{(n)})^{\alpha}+1]\\
&=& \frac{C(c_{\ref{def:singulartildek}}, L, \Psi, \Phi)  c_{0}(K)^{2\gamma}t}{n}
\frac{a_{n-1}^{(2\gamma-\frac{\alpha}{\xi} )} }{a_{n}}.
\end{eqnarray*}

\noindent
Observe now that
$\int_{a_{n}}^{a_{n-1}} x^{-1} dx \sim n$
so that $\frac{a_{n-1}}{a_{n}} \sim e^{n}$ or (using that $a_{0}=1$)
$a_{n} \sim e^{-\frac{n(n+1)}{2}}.$ Thus,
\begin{equation}
\label{eq:limI_3}
\lim_{n \rightarrow \infty} \IE \Big( |I_{3}^{m^{(n)},n}(t\wedge T_{\xi,K})| \Big)=0
\end{equation}
if
$n(n+1) - (2\gamma-\frac{\alpha}{\xi} )(n-1)n <0$
for $n$ large. This is equivalent to
\begin{equation*}
1-(2\gamma-\frac{\alpha}{\xi} ) <0 \Leftrightarrow
\alpha < \xi (2\gamma -1)
\end{equation*}
which holds by (\ref{alphbnd}).

Use (\ref{absconv}) and
Fatou's Lemma on the left-hand side of
(\ref{eq:Iparts}),
and Lemma~\ref{lem:conv} and (\ref{eq:limI_3}) on the
right-hand side, to take limits in this equation and
so conclude
\begin{eqnarray*}
\int_{\IR^{d}} \IE\Big( |\tilde{u}(t\wedge T_{\xi,K},x)|\Big) \Psi_t(x) dx
&\leq&  \liminf_{n \rightarrow \infty} \int_{\IR^{d}} \IE \Big(
\phi_{n}(\la \tilde{u}_{t \wedge T_{\xi,K}}, \Phi^{m^{(n)}}_{x} \ra) \Big)
\Psi_t(x) dx \\
&\leq&  \IE \Big( \int_{0}^{t \wedge T_{\xi,K}} \int_{\IR^{d}}
|\tilde{u}(s,x)|
\frac{1}{2}\Bigl(\Delta \Psi_s(x)+\dot\Psi_s(x)\Bigr)dx ds \Big)\\
&\leq& \int_{0}^{t} \int_{\IR^{d}} \IE \Big(  |\tilde{u}(s,x)| \Big)
|\frac{1}{2}\Delta \Psi_s(x)+\dot\Psi_s(x)|dx ds.
\end{eqnarray*}
Since $T_{\xi,K}$ tends in probability to infinity as $K \rightarrow \infty$
according to (\ref{eq:stoptbnd}), we know that \hfil\break
$\tilde{u}(t\wedge
T_{\xi,K},x)
\rightarrow \tilde{u}(t,x)$ and so we finally conclude with another application of
Fatou's Lemma that
\begin{equation*}
\int_{\IR^{d}} \IE\Big( |\tilde{u}(t,x)|\Big) \Psi_t(x) dx
\leq \int_{0}^{t} \int_{\IR^{d}} \IE \Big(  |\tilde{u}(s,x)| \Big)
\Bigl|\frac{1}{2}
\Delta \Psi_s(x)+\dot\Psi_s(x)\Bigr|dx ds.
\end{equation*}
This is (\ref{eq:heatGron}) of Section 3 and the conclusion now follows as in the proof of
Theorem~\ref{thm:uniquekbounded} given there.

\section{Proof of Theorem \ref{thm:smallHoelder}}
\label{proof:smallHoelder}

In this section we will first prove three technical lemmas
needed for the proof of Theorem \ref{thm:smallHoelder}.

\begin{lemma}
\label{lemma:correst}
Let $B$ be a standard d-dimensional Brownian motion.
For $\alpha<d$ there exists a constant
$c_{\ref{lemma:correst}}=c_{\ref{lemma:correst}}(\alpha,d)$ such that for 
all $x,y \in \IR^d$ and $t,t' >0$,
\begin{equation}
\label{eq:correst}
\int_{\IR^{d}}\int_{\IR^{d}}p_{t}(x-w) p_{t'}(y-z) |w-z|^{-\alpha} dw dz
= \IE_{x-y}(|B_{t+t'}|^{-\alpha}) \leq \IE_{0}(|B_{t+t'}|^{-\alpha})
\leq c_{\ref{lemma:correst}} (t+t')^{-\frac{\alpha}{2}}.
\end{equation}
 In addition, for
any $\lambda'>0$, $c\ge 0$, and $0<t\le t'$,
\begin{equation}
\label{eq:correstlambda}
\int_{\IR^{d}}\int_{\IR^{d}}e^{\lambda' (|w|+|z|)}
p_{t}(x-w) p_{t'}(y-z)[|w-z|^{-\alpha}+c] dw dz
\leq c_{\ref{lemma:correst}}e^{2(\lambda')^2t'} e^{\lambda' (|x|+|y|)}
[(t+t')^{-\frac{\alpha}{2}}+c].
\end{equation}
\end{lemma}
\begin{proof}
The first equality of (\ref{eq:correst}) is immediate from change of variables. The second
inequality then follows from a simple coupling argument: Let $|B^{i}_{t}|$ for
$i=1,2$ be the radial part of a d-dimensional Brownian motion started
at $0$ and $|x-y|$ respectively. Define the stopping time
$T:=\inf\{t \geq 0 : |B^1_t| > |B^2_t| \}.$ Then
\begin{equation*}
|B^3_t| =\left\{
\begin{array}{ll}
|B^2_t| & \text{ for } t \leq T,\\
|B^1_t| & \text{ for } t > T,
\end{array}\right.
\end{equation*}
has the same law as $|B^2|$ and the property that $|B^3_{t}|\geq
|B^1_{t}|$ for all $t\geq 0$  a.s., which implies the inequality
of the expectations in (\ref{eq:correst}). We finally compute
by setting $r=\frac{|w|^2}{2t},$
\begin{eqnarray*}
\IE_{0}(|B_{t}|^{-\alpha})
= \int_{\IR^{d}}|w|^{-\alpha} (2 \pi t)^{-\frac{d}{2}} \exp(-\frac{|w|^{2}}{2t}) dw
= c_d \int_{0}^{\infty} r^{\frac{d-\alpha}{2}-1}
\exp(-r) dr  \cdot t^{-\frac{\alpha}{2}}
= c(\alpha,d) t^{-\frac{\alpha}{2}}
\end{eqnarray*}
provided that $\alpha<d.$ This shows (\ref{eq:correst}). For proving (\ref{eq:correstlambda}) we
note that for $0<t\le t',$
\begin{equation}
\label{eq:expest}
e^{\lambda' |w|} p_{t}(w) \leq 2^{\frac{d}{2}} \exp( \lambda' |w| -\frac{1}{4t} |w|^{2})
p_{2t}(w) \leq c_de^{(\lambda')^2t'}p_{2t}(w)
\end{equation}
since $\lambda' |w| -\frac{1}{4t} |w|^{2} \leq (\lambda')^2 t.$  Therefore,
\begin{eqnarray*}
& &\int_{\IR^{d}}\int_{\IR^{d}}e^{\lambda' (|w|+|z|)} p_{t}(x-w) p_{t'}(y-z)
[ |w-z|^{-\alpha}+c] dw dz\\
&\leq&  e^{\lambda' (|x|+|y|)}
\int_{\IR^{d}}\int_{\IR^{d}}e^{\lambda' (|w|+|z|)} p_{t}(w) p_{t'}(z) 
[|w-z+x-y|^{-\alpha}+c]
dw dz\\ &\leq&  c_d^2e^{2(\lambda')^2t'}e^{\lambda' (|x|+|y|)}
\int_{\IR^{d}}\int_{\IR^{d}}p_{2t}(w) p_{2t'}(z) [|w-z+x-y|^{-\alpha}+c] dw dz\\
&\leq& c(\alpha,d)e^{2\lambda^2t'} e^{\lambda (|x|+|y|)}[  (t+t')^{-\frac{\alpha}{2}}+c].
\end{eqnarray*}
Here, we have used a shift of variables in the first and 
(\ref{eq:expest}) in the second
inequality as well as (\ref{eq:correst}) in the third. This shows (\ref{eq:correstlambda}).
\end{proof}

The next lemma provides some estimates of the temporal and spatial
differences of the heat kernels:
\begin{lemma}
\label{lemma:pdiffest}
There are constants  $c_{\ref{eq:pdiffest}}(d)$ and $c_{\ref{eq:spacecorrest}}(\alpha,d)$
such that if $0 < \beta \leq 1$ and 
$\lambda' \geq 0$, then for any $x,y \in \IR^d,$ $0<t\leq
t',$
\begin{eqnarray}
\label{eq:pdiffest}
\int_{\IR^{d}}|p_{t}(x-w) -p_{t'}(y-w)| e^{\lambda' |w|} dw
&\leq&  c_{\ref{eq:pdiffest}}e^{2(\lambda')^2t'}
\left(e^{\lambda' |x|} + e^{\lambda' |y|} \right)
e^{2\beta\lambda'(|x-y|)} \\
\nonumber
& &\phantom{AAAAAA}
\left( t^{-\beta/2} |x-y|^{\beta}
+ t^{-\beta} |t'-t|^{\beta}  \right),
\end{eqnarray}
\begin{eqnarray}
\label{eq:spacecorrest}
\nonumber
&&\int_{\IR^{d}}\int_{\IR^{d}}|p_{t}(x-w) -p_{t}(y-w)|
|p_{t}(x-z) -p_{t}(y-z)| [|w-z|^{-\alpha} +1]\,dw dz\\
&&\qquad\leq c_{\ref{eq:spacecorrest}} [t^{-1-\frac{\alpha}{2}}+t^{-1}] |x-y|^{2},
\end{eqnarray}
and
\begin{eqnarray}
\label{eq:timecorrest}
\nonumber
&&\int_{\IR^{d}}\int_{\IR^{d}}|p_{t}(x-w) -p_{t'}(x-w)|
|p_{t}(x-z) -p_{t'}(x-z)| [|w-z|^{-\alpha}+1]\, dw dz\\
&&\qquad\leq c_{\ref{eq:spacecorrest}} [t^{-2-\frac{\alpha}{2}} +t^{-2}]|t'-t|^{2}.
\end{eqnarray}
\end{lemma}
\begin{proof}
We consider the space and time differences separately.
For the former, define $v=x-y,$ and set $\hat{v}_{0}=0, \hat{v}_{d}=v,$
and $\hat{v}_{i}-\hat{v}_{i-1}=v_{i} e_{i},$ where $v_{i}$ is the i-th component
of $v$ and $e_{i}$ is the i-th unit vector in $\IR^{d}.$ Therefore,
\begin{eqnarray*}
|\exp(-\frac{|w+v|^{2}}{2t}) - \exp(-\frac{|w|^{2}}{2t})|
&\leq&\sum_{i=1}^{d}
|\exp(-\frac{|w+\hat{v}_{i}|^{2}}{2t}) - \exp(-\frac{|w+\hat{v}_{i-1}|^{2}}{2t})|\\
&=& \sum_{i=1}^{d} |\int_{0}^{v_{i}} \frac{w_{i}+r_{i}}{t}
\exp(-\frac{|w+\hat{v}_{i-1}+r_{i}e_i|^{2}}{2t}) dr_{i}|.
\end{eqnarray*}
Hence, by a change of variables, (\ref{eq:expest}), and using $|w|\le |\hat w_i|+|w_i|$
($\hat w_i=w-w_ie_i$), we have
\begin{eqnarray}
\nonumber
& &\int_{\IR^{d}}|p_{t}(x-w) -p_{t}(y-w)| e^{\lambda' |w|} dw \\
\nonumber
&\leq& e^{\lambda'|x|} (2 \pi t)^{-\frac{d}{2}} \sum_{i=1}^{d} \int_{0}^{|v_{i}|}
\int_{\IR^{d}}\frac{|w_{i}+r_{i}|}{t}
\exp(-\frac{|w+\hat{v}_{i-1}+r_{i}e_i|^{2}}{2t})e^{\lambda' |w|}dw dr_i\\
\nonumber
&\leq & c_de^{(\lambda')^2t'} e^{\lambda'(|x|+|v|)}
t^{-\frac{1}{2}} \sum_{i=1}^{d} \int_{0}^{|v_{i}|}
\int_{-\infty}^{\infty} \frac{|w_{i}+r_{i}|}{t}
\exp(-\frac{(w_{i}+r_{i})^{2}}{2t}) e^{\lambda'  |w_{i}|} dw_{i}  dr_{i}\\
\nonumber
&\leq&  c_de^{(\lambda')^2t'}e^{\lambda'(|x|+|v|)}
t^{-\frac{1}{2}} \left( \sum_{i=1}^{d} e^{\lambda' |v_i|}|v_{i}| \right)
\int_{0}^{\infty} \frac{r}{t} \exp(-\frac{r^{2}}{4t}) dr\\
\label{eq:spacediff}
&\leq&  c_de^{(\lambda')^2 t'}e^{\lambda'|x|+2\lambda'|x-y|} t^{-\frac{1}{2}} |x-y|.
\end{eqnarray}
Similarly, using that $a \leq c\exp(a^2/4)$ for all $a \in \IR_+,$ we get
\begin{eqnarray*}
& &\int_{\IR^{d}}\int_{\IR^{d}}|p_{t}(x-w) -p_{t}(y-w)|
\cdot|p_{t}(x-z) -p_{t}(y-z)| [|w-z|^{-\alpha}+1]\, dw dz \\
&\leq&  C t^{-d}\int_{\IR^{d}}\int_{\IR^{d}}
\sum_{i,j=1}^{d} \int_{0}^{|v_{i}|} \int_{0}^{|v_{j}|} \frac{|w_{i}+r_{i}|}{t}
\exp(-\frac{|w+\hat{v}_{i-1}+r_{i}e_i|^{2}}{2t}) \\
& & \phantom{AAAAAAAAAAAAAA}
\frac{|z_{i}+\tilde{r}_{i}|}{t}
\exp(-\frac{|z+\hat{v}_{i-1}+\tilde{r}_{i}e_i|^{2}}{2t})  dr_{i} d\tilde{r}_{i}
[|w-z|^{-\alpha}+1] \,dw dz\\
&\leq&  C t^{-d-1}
\sum_{i,j=1}^{d}  \int_{0}^{|v_{i}|} \int_{0}^{|v_{j}|}
\big( \int_{\IR^{d}}\int_{\IR^{d}}
\exp(-\frac{|w+\hat{v}_{i-1}+r_{i}e_i|^{2}}{2t}+\frac{|w_{i}+r_{i}|^{2}}{4t}) \\
& & \phantom{AAAAAAAAAAAAAA}
\exp(-\frac{|z+\hat{v}_{i-1}+\tilde{r}_{i}e_i|^{2}}{2t}+ \frac{|z_{i}+\tilde{r}_{i}|^2}{4t})
 [|w-z|^{-\alpha} +1]\,dw dz \big)dr_{i} d\tilde{r}_{i}\\
&\leq&  C(\alpha,d) t^{-1}
\sum_{i,j=1}^{d}  \int_{0}^{|v_{i}|} \int_{0}^{|v_{j}|}
[t^{-\frac{\alpha}{2}} +1]dr_{i} d\tilde{r}_{i}\\
&\leq& C(\alpha,d) [t^{-1-\frac{\alpha}{2}}+t^{-1}] |x-y|^{2},
\end{eqnarray*}
where we have used an appropriate shift of variables and Lemma \ref{lemma:correst}
in the previous to last line. This shows (\ref{eq:spacecorrest}).

For the time differences observe that for some $C=C(\alpha,d)$,
\begin{eqnarray}
\nonumber
|p_{t}(w) -p_{t'}(w)| &\leq& C |t^{-\frac{d}{2}} -
{t'}^{-\frac{d}{2}}| \exp(-\frac{|w|^{2}}{2t}) + C
{t'}^{-\frac{d}{2}} |\exp(-\frac{|w|^{2}}{2t})-
\exp(-\frac{|w|^{2}}{2t'})|\\
\nonumber
&\leq& C |t'-t| t^{-\frac{d}{2}-1} \exp(-\frac{|w|^{2}}{2t}) + C
{t'}^{-\frac{d}{2}} \int_{t}^{t'}
\exp(-\frac{|w|^{2}}{2s}) \frac{|w|^{2}}{2s^2} ds \\
\label{eq:ptdiff}
&\leq& C t^{-1} |t'-t|     \left( p_{t}(w) + p_{2t'}(w) \right),
\end{eqnarray}
since $\frac{|w|^2}{4s}\leq \exp(\frac{|w|^2}{4s}).$
Therefore, another application of (\ref{eq:expest}) yields
\[\int_{\IR^{d}}|p_{t}(w) -p_{t'}(w)| e^{\lambda' |w|} dw 
\leq  C(d) e^{(\lambda')^2t'} t^{-1} |t'-t|.\]
Taking this estimate with a change of variables, together with (\ref{eq:spacediff}),
we obtain
\begin{equation}\int|p_t(x-w)-p_{t'}(y-w)|e^{\lambda'|w|}\,dw\le
C(d) e^{2\lambda^2t'}e^{\lambda'|x|}\Bigl[e^{2\lambda'|x-y|}{|x-y|\over \sqrt t}
+{|t'-t|\over t}\Bigr].\end{equation}
An application of (\ref{eq:expest}) and a change of variables also shows that 
\begin{equation}
\int|p_t(x-w)-p_{t'}(y-w)|e^{\lambda'|w|}dw\le
C(d) e^{(\lambda')^2t'}(e^{\lambda'|x|}+e^{\lambda'|y|}).
\end{equation}
If $\beta\in(0,1]$, the inequality $z\wedge 1\le z^\beta$ for $z\ge 0$, and the previous
two bounds now show that
\begin{eqnarray*}
\int|p_t(x-w)-p_{t'}(y-w)|e^{\lambda'|w|}dw&\le&C(d) e^{2(\lambda')^2t'}
(e^{\lambda'|x|}+e^{\lambda'|y|})\\
& &\quad\times[e^{2\beta\lambda'|x-y|}|x-y|^\beta t^{-\beta/2}+|t'-t|^\beta
t^{-\beta}],
\end{eqnarray*}
which implies (\ref{eq:pdiffest}).
Similarly, using (\ref{eq:ptdiff}) and (\ref{eq:correst}),
\begin{eqnarray*}
& &\int_{\IR^{d}}\int_{\IR^{d}}|p_{t}(x-w) -p_{t'}(x-w)|
\cdot|p_{t}(x-z) -p_{t'}(x-z)| [|w-z|^{-\alpha} +1]\,dw dz \\
&\leq&  C t^{-2} |t'-t|^{2} \int_{\IR^{d}}\int_{\IR^{d}}(p_{t}(w) + p_{2t'}(w))
(p_{t}(z) + p_{2t'}(z))
[|w-z|^{-\alpha}+1]\, dw dz\\
&\leq&  C(\alpha,d)( t^{-2-\frac{\alpha}{2}}+t^{-2}) |t'-t|^{2},
\end{eqnarray*}
which proves (\ref{eq:timecorrest}).
\end{proof}

We will also need the following rather technical lemma:

\begin{lemma}
\label{lemma:Jest}
For 
$b,c\geq 0$ with  $c< \frac{1}{2}(b+ 1-\frac{\alpha}{2})$, and $a\in(c,1-\alpha/2)$,
there is a finite constant $c_{\ref{lemma:Jest}}=c_{\ref{lemma:Jest}}(a,b,c,\alpha)$ such
that 
$t\geq 0$,
\begin{eqnarray*}
Q(t,a,b,c,\alpha)&:=&\int_{0}^{t} \int_{0}^{t}
 (t-r)^{a-1-c}  (t-r')^{a-1-c}
 \int_{0}^{r \wedge r'}   (t-s)^{b} (r-s)^{-a} (r'-s)^{-a}\\
& &\phantom{AAA}
 \int_{\IR^{d}}\int_{\IR^d} p_{r-s}(w) p_{r'-s}(z) [|w-z|^{-\alpha}+1]\,
 dw dz ds dr dr' \\
&\le &c_{\ref{lemma:Jest}}[t^{b+1-\alpha/2-2c}+t^{b+1-2c}].
\end{eqnarray*}
\end{lemma}
\begin{proof}
By Lemma \ref{lemma:correst} it suffices to estimate
\begin{eqnarray*}
& &\int_{0}^{t} \int_{0}^{t}
(t-r)^{a-1-c} (t-r')^{a-1-c}\\
& &\phantom{AAAAAAAAA}
\int_{0}^{r \wedge r'} (t-s)^{b} (r-s)^{-a} (r'-s)^{-a}
[(r-s+r'-s)^{-\frac{\alpha}{2}}+1]\,
ds dr dr'\\
&=& 2 \int_{0}^{t} \int_{s}^{t} (t-r)^{a-1-c} (t-s)^{b} (r-s)^{-a}\\
&&\phantom{AAAAAAAAAAAA}\times
\left(\int_{r}^{t} (t-r')^{a-1-c}  [(r'-s)^{-a-\frac{\alpha}{2}}+(r'-s)^{-a}]\, dr'\right)
dr ds ,
\end{eqnarray*}
where we have used the symmetry in $r$ and $r'$ and concentrated on the case $r \leq r'.$
Substituting $v=\frac{r'-r}{t-r}$ and using that $c<a < 1-\frac{\alpha}{2}$ we calculate
for $t\ge r\ge s$,
\begin{eqnarray*}
&&\int_{r}^{t} (t-r')^{a-1-c}  [(r'-s)^{-a-\frac{\alpha}{2}}+(r'-s)^{-a}]\, dr'\\
&&\quad= (t-r)^{-\frac{\alpha}{2}-c} \int_{0}^{1}
(1-v)^{a-1-c} \Bigl[(v+\frac{r-s}{t-r})^{-a-\frac{\alpha}{2}}
+(t-r)^{\frac{\alpha}{2}}(v+\frac{r-s}{t-r})^{-a}\Bigr]\, dv
\\ 
&&\quad\leq  C(a,c,\alpha) 
(t-r)^{a-c}[(t-r)^{-a-\alpha/2}\wedge(r-s)^{-a-\alpha/2}+(t-r)^{-a}\wedge(r-s)^{-a}].
\end{eqnarray*}
Hence the required $Q$ is at most $C(a,c,\alpha)$ times the sum of the following integral,
$I(\beta)$, for $\beta=a$ and $\beta=a+\alpha/2$:
\[I(\beta)=\int_{0}^{t} (t-s)^{b} \int_{s}^{t}
(t-r)^{2a-1-2c} (r-s)^{-a}
\left( (t-r)^{-\beta} \wedge (r-s)^{-\beta} \right)
 dr ds.
\]
For these values of $\beta$, $I(\beta)$ is at most
\begin{eqnarray*}
& &\int_{0}^{t} (t-s)^{b}\left( \int_{s}^{\frac{t+s}{2}}
(t-r)^{2a-1-2c-\beta} (r-s)^{-a} dr + \int_{\frac{t+s}{2}}^{t}
(t-r)^{2a-1-2c} (r-s)^{-a-\beta} dr  \right)  ds\\
&\leq& C(a,\alpha) \int_{0}^{t} \left(
(t-s)^{b+2a-1-2c-\beta}\int_{s}^{\frac{t+s}{2}} (r-s)^{-a} dr +
(t-s)^{b-a-\beta} \int_{\frac{t+s}{2}}^{t} (t-r)^{2a-1-2c} dr  \right)  ds\\
&\leq& C(a,c,\alpha) \int_{0}^{t}  (t-s)^{b+a-2c-\beta} ds \leq
C(a,b,c,\alpha)t^{b+a-2c-\beta+1}.
\end{eqnarray*}
Here, we have used that $t-r\geq \frac{t-s}{2}$ for $r \in [s,\frac{t+s}{2}]$ and analogously
$r-s\geq \frac{t-s}{2}$ for $r \in [\frac{t+s}{2}, t]$ as well as our assumption of $a>c$ and
$c< \frac{1}{2}(b+ 1-\frac{\alpha}{2}).$  The result follows upon summing over the two
values of $\beta$.
\end{proof}

\noindent
PROOF OF THEOREM \ref{thm:smallHoelder}.

\noindent
Fix arbitrary (deterministic) $(t,x), (t',y)$ such that
$d((t,x),(t',y))\le \epsilon\equiv 2^{-N}$ ($N\in\IN$) and $t\le t'$ (the case $t' \leq t$
works analogously).
As $\xi_1<(\xi\gamma+1-\alpha/2)\wedge1$, we may choose
$\delta\in(0,1-\alpha/2)$ so that
\begin{equation}
\label{chdelta}
1>\xi\gamma+1-\alpha/2-\delta>\xi_1.
\end{equation}
Note that $\xi\gamma<1$ shows we may choose $\delta$ in the required  range. Next choose
$\delta'\in(0,\delta)$ and $p\in(0,\xi\gamma)$ so that
\begin{equation}
\label{chp}
1>p+1-\alpha/2-\delta>\xi_1,
\end{equation}
and
\begin{equation}
\label{chdelta'}
1>\xi\gamma+1-\alpha/2-\delta'>\xi_1.
\end{equation}

Now consider for some random $N_1=N_1(\omega,\xi,\xi_1)$ to be chosen below,
\begin{eqnarray}
\label{eq:H1}
& &\IP\left( |\tilde u(t,x)-\tilde u(t,y)| \geq
|x-y|^{1-\frac{\alpha}{2} -\delta}\epsilon^{p},  (t,x)\in Z_{K,N,\xi}, N\ge N_1\right)\\
\nonumber
&+&\IP\left( |\tilde u(t',x)-\tilde u(t,x)| \geq
|t'-t|^{\frac{1}{2}(1-\frac{\alpha}{2}-\delta)}\epsilon^{p}, (t,x)\in Z_{K,N,\xi}, t'\le T_K,
N\ge N_1
 \right).
\end{eqnarray}
In order to simplify notation we define
\begin{eqnarray*}
D^{x,y,t,t'}(w,z,s)&=&
\left| p_{t-s}(x-w)- p_{t'-s}(y-w)\right|
\left| p_{t-s}(x-z)- p_{t'-s}(y-z)\right|\\
& &
\cdot |\tilde u(s,w)|^{\gamma}|\tilde u(s,z)|^{\gamma} [|w-z|^{-\alpha}+1],\\
D^{x,t'}(w,z,s)&=& p_{t'-s}(x-w) p_{t'-s}(x-z)
|\tilde u(s,w)|^{\gamma}|\tilde u(s,z)|^{\gamma} [|w-z|^{-\alpha}+1].
\end{eqnarray*}
With this notation expression (\ref{eq:H1}) is bounded by
\begin{eqnarray}
\label{eq:H3}
& &\IP\Big( |\tilde{u}(t,x)-\tilde{u}(t,y)|
\geq |x-y|^{1-\frac{\alpha}{2} -\delta}\epsilon^{p}, (t,x)\in Z_{K,N,\xi},N \geq N_1\\
\nonumber
& &\phantom{AAAAAAAAAAAAAAAAAAAAAA}
\int_{0}^{t} \int_{\IR^{d}}\int_{\IR^{d}} D^{x,y,t,t}(w,z,s)  dw dz ds
\leq  |x-y|^{2-\alpha-2\delta'}\epsilon^{2p} \Big)\\
\nonumber
&+&\IP\Big( |\tilde{u}(t',x)-\tilde{u}(t,x)|
\geq |t'-t|^{\frac{1}{2}(1-\frac{\alpha}{2}-\delta)}\epsilon^{p}, (t,x)\in Z_{K,N,\xi}, t'\le
T_K, N \geq N_1\\
\nonumber
& &\phantom{A}
\int_{t}^{t'} \int_{\IR^{d}}\int_{\IR^{d}} D^{x,t'}(w,z,s) dw dz ds
+ \int_{0}^{t} \int_{\IR^{d}}\int_{\IR^{d}} D^{x,x,t,t'}(w,z,s)
dw dz ds \leq
(t'-t)^{1-\frac{\alpha}{2}-\delta'}\epsilon^{2p}
\Big)\\
\nonumber
&+& \IP\Big( \int_{0}^{t} \int_{\IR^{d}}\int_{\IR^{d}}
D^{x,y,t,t}(w,z,s)
dw dz ds
> |x-y|^{2-\alpha-2\delta'}\epsilon^{2p}, (t,x)\in
Z_{K,N,\xi}, N\ge N_1
\Big)\\
\nonumber
&+& \IP\Big(  \int_{t}^{t'} \int_{\IR^{d}}\int_{\IR^{d}}
D^{x,t'}(w,z,s) dw dz ds
+ \int_{0}^{t} \int_{\IR^{d}}\int_{\IR^{d}} D^{x,x,t,t'}(w,z,s) dw dz ds\\
\nonumber
& &\phantom{AAAAAAAAAAAAAAAAAAAA}
>  (t'-t)^{1-\frac{\alpha}{2}-\delta'}\epsilon^{2p},
 (t,x)\in Z_{K,N,\xi}, t'\le T_K, N\ge N_1 \Big)\\
\nonumber
&=:& P_{1}+P_{2}+P_{3}+P_{4}.
\end{eqnarray}
Notice that the processes $\tilde{t} \mapsto \int_{0}^{\tilde{t}}\int_{\IR^{d}} p_{t-s}(x-w)
\left(\sigma(u^{1}(s,w))-\sigma(u^{2}(s,w)\right) W(dw ds)$
are continuous local martingales for any fixed $x,t$ on $0 \leq \tilde{t} \leq t$. 
We bound the appropriate differences of these integrals by considering the respective
quadratic variations of $\tilde{u}(t,x)-\tilde{u}(t,y)$  and
$\tilde{u}(t',x)-\tilde{u}(t,x)$ (see (\ref{heatecolint})).  If
$|\sigma(u)-\sigma(v)|\le L|u-v|^\gamma$ and recalling that $|k(x,y)|\le
c_{\ref{thm:uniquekunbounded}}[|x-y|^{-\alpha}+1]$, we see that the time integrals in the above
probabilities differ from the appropriate square functions by a multiplicative factor of
$L^2 c_{\ref{thm:uniquekunbounded}}$. 

If $\delta''=\delta-\delta'>0$, $B$ is a standard one-dimensional Brownian motion with
$B(0)=0$, and $B^{*}(t):=\sup_{0\leq s \leq t} |B(s)|,$ then
the first two probabilities of (\ref{eq:H3}) can be bounded
using the Dubins-Schwarz Theorem:
\begin{eqnarray}
\label{P1est}
\nonumber
P_{1}&\leq& \IP\left( B^{*}(c_{\ref{thm:uniquekunbounded}}
L^2|x-y|^{2-\alpha-2\delta'}\epsilon^{2p})
\geq |x-y|^{1-\frac{\alpha}{2} -\delta}\epsilon^{p}\right)\\
\nonumber
&=& \IP\left( B^{*}(1) \sqrt{c_{\ref{thm:uniquekunbounded}}} L
|x-y|^{1-\frac{\alpha}{2}-\delta'}\epsilon^{p}
\geq |x-y|^{1-\frac{\alpha}{2} -\delta}\epsilon^{p}
\right)\\
&=& \IP\left( B^{*}(1)
\geq (\sqrt{c_{\ref{thm:uniquekunbounded}}} L)^{-1} |x-y|^{-\delta'' } \right)
\leq c_{\ref{P1est}} \exp(-c'_{\ref{P1est}} |x-y|^{-\delta''}),
\end{eqnarray}
where we have used the reflection principle in the last line.
Likewise,
\begin{eqnarray}
\label{P2est}
\nonumber
P_{2}&\leq& \IP\left(
B^{*}(c_{\ref{thm:uniquekunbounded}}L^2|t'-t|^{1-\frac{\alpha}{2}-\delta'}\epsilon^{2p})
\geq |t'-t|^{\frac{1}{2}(1-\frac{\alpha}{2}-\delta)}\epsilon^{p}\right)\\
&=& \IP\left( B^{*}(1)
\geq  (\sqrt{c_{\ref{thm:uniquekunbounded}}} L)^{-1}|t'-t|^{-\frac{\delta''}{2} } \right)
\leq c_{\ref{P1est}} \exp(-c'_{\ref{P1est}} |t'-t|^{-\frac{\delta''}{2}}).
\end{eqnarray}
Here the constants $c_{\ref{P1est}}$  and
$c'_{\ref{P1est}}$ depend on $d$, $L$, and $c_{\ref{thm:uniquekunbounded}}$. 

In order to bound $P_{3}$ and $P_{4}$ we
estimate the respective integral expressions by splitting them up in several parts: Let
$\delta_1\in (0,\frac{1}{2}(1-\frac{\alpha}{2}))$ and $t_{0}=0, t_{1}=t-\epsilon^2,
t_{2}=t$ and
$t_3=t'.$ We also define
\begin{eqnarray}
A_{1}^{1,s}(x)
&=&\{ w \in \IR^{d} : |x-w| \leq 2 \sqrt{t-s} \epsilon^{-\delta_{1}}\} \hbox{ and }
A_{2}^{1,s}(x)= \IR^{d}\setminus A_{1}^{1,s}(x),\\
A_{1}^{2}(x) &=& \{ w \in \IR^{d} : |x-w| \leq 2 \epsilon^{1-\delta_{1}}\} \hbox{ and }
A_{2}^{2}(x)= \IR^{d}\setminus A_{1}^{2}(x).
\end{eqnarray}
For notational convenience we will sometimes omit the index $s$ for $A^{1}_{i}(x)$. We continue to
write
\begin{equation*}
Q^{x,y,t,t'}:= \int_{0}^{t} \int_{\IR^{d}}\int_{\IR^{d}} D^{x,y,t,t'}(w,z,s)
dw dz ds = \sum_{i,j,k=1,2} Q_{i,j,k}^{x,y,t,t'},
\end{equation*}
where
\begin{equation*}
Q_{i,j,k}^{x,y,t,t'}
:=\int_{t_{i-1}}^{t_{i}} \int_{A^{i}_{j}(x)}\int_{A^{i}_{k}(x)} D^{x,y,t,t'}(w,z,s)
dw dz ds.
\end{equation*}
And likewise,
\begin{equation*}
Q^{x,t,t'}:= \int_{t}^{t'} \int_{\IR^{d}}\int_{\IR^{d}} D^{x,t'}(w,z,s)
dw dz ds = \sum_{j,k=1,2} Q_{j,k}^{x,t,t'},
\end{equation*}
where
\begin{equation*}
Q^{x,t,t'}_{j,k}
:=\int_{t}^{t'} \int_{A^{2}_{j}(x)}\int_{A^{2}_{k}(x)} D^{x,t'}(w,z,s)
dw dz ds.
\end{equation*}
Before we proceed let us note that $\tilde{u}$ can be bounded on the sets $A^{i}_{1}$
as follows: Set
\begin{eqnarray}
\label{N1bnd}
N_1(\omega)=\Bigl[{5N_\xi(\omega)\over \delta_1}\Bigr]\ge
\Bigl[{N_\xi(\omega)+4\over 1-\delta_1}\Bigr]\in\IN,
\end{eqnarray}
where $[\cdot]$ is the
greatest integer function and assume $N\ge N_1$ in the following.

Recall $\lambda>0$ is a fixed constant used in the definition of $T_K$ and hence
$Z_{K,N,\xi}$.  As it is fixed, we often suppress dependence on $\lambda$ in our notation.
\begin{lemma}
\label{lemma:tubnd}
Let $N\ge N_1.$ Then on $\{\omega:(t,x)\in Z_{K,N,\xi}\}$,
\begin{eqnarray}
\label{tubnd}
|\tilde u(s,w)| &\le& 10\epsilon^{(1-\delta_1)\xi}
\quad \quad \quad \text{ for } s \in [t-\epsilon^{2},t'], w \in A_{1}^{2}(x),\\
\label{tubnd4}
|\tilde u(s,w)| &\leq&  (8 + 3 K  2^{N_{\xi}\xi})
e^{\lambda |w|} (t-s)^{\frac{\xi}{2}}\epsilon^{-\delta_1 \xi}
\quad \quad \text{ for } s \in [0,t-\epsilon^2], w \in A_{1}^{1,s}(x).
\end{eqnarray}
\end{lemma}
\begin{proof}
We choose $N'\in\IN$ so that
$2^{-N'-1}\le 3\epsilon^{1-\delta_1}\le 2^{-N'}$.  Then
$2^{-N'-3}<2^{-N(1-\delta_1)}\le 2^{-N'-1},$
and so by (\ref{N1bnd}),
\begin{eqnarray}
\label{N'bnd}
N' > N(1-\delta_1)-3\ge N_1(1-\delta_1)-3\ge N_\xi.
\end{eqnarray}

\noindent
Assume $(t,x)\in Z_{K,N,\xi}$, $0\le t'\le T_K$ and choose $(\hat{t},\hat{x})$ such that
\begin{eqnarray}
\label{t0x0}
\hbox{
$\hat{t}\le T_K$, $d((t,x),(\hat{t},\hat{x}))<\epsilon=2^{-N}$, and $|\tilde u(\hat{t},\hat{x})|\le
2^{-N\xi}=\epsilon^\xi$.}
\end{eqnarray}
We first observe that for $s \in [t-\epsilon^{2},t']$
and $w \in A_{1}^{2}(x)$ so that $|w-x|\leq 2\epsilon^{1-\delta_1},$ we have
\begin{eqnarray}
\label{swtx}
d((s,w),(t,x))\le \epsilon+2\epsilon ^{1-\delta_1}\le
3\epsilon^{1-\delta_1}\le 2^{-N'}.
\end{eqnarray}
Therefore by (\ref{reghyp}) and (\ref{N'bnd}), for $s \in [t-\epsilon^{2},t']$
and $w \in A_{1}^{2}(x),$
\begin{eqnarray}
\nonumber
|\tilde u(s,w)|&\le & |\tilde u(\hat{t},\hat{x})|+|\tilde u(\hat{t},\hat{x})-\tilde u(t,x)|+|\tilde
u(t,x)-\tilde u(s,w)|\\
\nonumber
&\le &  2\cdot2^{-N\xi}+2^{-N'\xi}\\
\nonumber
&\le &2\epsilon^\xi +(8\epsilon^{1-\delta_1})^\xi\\
\label{tubnd1}
&\le &10\epsilon^{(1-\delta_1)\xi},
\end{eqnarray}
which proves (\ref{tubnd}).
Similarly, if $s \in [0,t-\epsilon^2]$ and $w \in A_{1}^{1,s}(x)$ meaning
that $|w-x|\leq  2\sqrt{t-s} \epsilon^{-\delta_1},$ we have
\begin{eqnarray}
\label{swtx2}
d((s,w),(t,x))\le \sqrt{t-s} +2\sqrt{t-s} \epsilon ^{-\delta_1}\le
3\sqrt{t-s}\epsilon^{-\delta_1}.
\end{eqnarray}
Notice that if $3\sqrt{t-s}\epsilon^{-\delta_1} \leq 2^{-N_{\xi}}$ then there exists an
$N' \geq N_{\xi}$  such that $  2^{-(N'+1)}\leq 3\sqrt{t-s}\epsilon^{-\delta_1}
\leq 2^{-N'}$  so that we can as in (\ref{tubnd}) bound
\begin{eqnarray}
\nonumber
|\tilde u(s,w)|&\le & |\tilde u(\hat{t},\hat{x})|+|\tilde u(\hat{t},\hat{x})-\tilde u(t,x)|
+|\tilde u(t,x)-\tilde{u}(s,w)|\\
\nonumber
&\leq& 2^{-N\xi} +2^{-N\xi} + 2^{-N'\xi}\\
\nonumber
&\leq &  2\cdot2^{-N\xi}+2\cdot 2^{-(N'+1)\xi}\\
\nonumber
&\leq &2 (t-s)^{\frac{\xi}{2}}
+ 2 \cdot 3^{\xi} (t-s)^{\frac{\xi}{2}} \epsilon^{-\delta_1 \xi} \\
\label{tubnd2}
&\leq & 8 (t-s)^{\frac{\xi}{2}} \epsilon^{-\delta_1 \xi},
\end{eqnarray}
since $\epsilon=2^{-N}\leq \sqrt{t-s}.$ If on the other hand
$3\sqrt{t-s}\epsilon^{-\delta_1} > 2^{-N_{\xi}}$ then we bound
\begin{eqnarray}
\nonumber
|\tilde u(s,w)| &\leq & K e^{\lambda |w|} \\
\nonumber
&=& (K  (t-s)^{-\frac{\xi}{2} }) e^{\lambda |w|} (t-s)^{\frac{\xi}{2}}\\
\label{tubnd3}
&\leq& (K  \epsilon^{-\delta_{1}\xi} 3^{\xi} 2^{N_{\xi}\xi})
e^{\lambda |w|} (t-s)^{\frac{\xi}{2}}.
\end{eqnarray}
Taking (\ref{tubnd2}) and (\ref{tubnd3}) together we obtain
(\ref{tubnd4}).
\end{proof}

\noindent
In the rest of this section $C(K)$ denotes a constant depending on $K$ (and possibly
$\lambda$) which may change from line to line.  We will first consider the terms for which
$j=k=1$ so that we can use the bounds (\ref{tubnd}) and (\ref{tubnd4}) of Lemma
\ref{lemma:tubnd}:
\begin{lemma}
\label{lemma:Q11}
If $0<\beta< 1-\frac{\alpha}{2}$,
$\beta' < \xi \gamma + 1-\frac{\alpha}{2},$ and $\beta'\le 1$, then on $\{\omega:(t,x)\in
Z_{K,N,\xi}\}$, 
\begin{eqnarray}
\label{Qx211}
Q_{2,1,1}^{x,y,t,t}
&\leq& c_{\ref{Qx211}}(\alpha,d,\beta,K)\epsilon^{2(1-\delta_1)\xi \gamma} |x-y|^{2\beta},\\
\label{Qt211}
Q_{2,1,1}^{x,x,t,t'}
&\leq& c_{\ref{Qx211}}(\alpha,d,\beta,K) \epsilon^{2(1-\delta_1)\xi \gamma}
|t'-t|^{\beta},\\
\label{Qx111}
Q_{1,1,1}^{x,y,t,t}
&\leq& c_{\ref{Qx111}}(\alpha,d,\beta',\xi\gamma,K) (8 + 3 K 2^{  N_{\xi} \xi})^{2\gamma} 
\epsilon^{-2\delta_1\xi \gamma }
 |x-y|^{2\beta'},\\
\label{Qt111}
Q_{1,1,1}^{x,x,t,t'}&\leq&
c_{\ref{Qx111}}(\alpha,d,\beta',\xi\gamma,K)(8 + 3 K 2^{  N_{\xi} \xi})^{2\gamma}  
\epsilon^{-2\delta_1\xi \gamma } |t'-t|^{\beta'},\\
\label{Qt11}
Q^{x,t,t'}_{1,1} &\leq& c_{\ref{Qt11}}(\alpha,d) \epsilon^{2\gamma \xi(1-\delta_{1})}
|t'-t|^{1-\frac{\alpha}{2}}.
\end{eqnarray}
\end{lemma}

\begin{proof}
Using the bounds (\ref{tubnd}) and (\ref{tubnd4}) of Lemma \ref{lemma:tubnd}
we obtain
\begin{eqnarray}
\label{eq:uest1}
Q_{2,1,1}^{x,y,t,t'} &\leq& 100^\gamma \epsilon^{2(1-\delta_1)\xi \gamma}
 \int_{t-\epsilon^{2}}^{t}
\int_{A^{2}_{1}(x)}\int_{A^{2}_{1}(x)} \left| p_{t-s}(x-w)- p_{t'-s}(y-w)\right|\\
\nonumber
&&\phantom{AAAAAAAAAAAAA}
\cdot
\left| p_{t-s}(x-z)- p_{t'-s}(y-z)\right|  [|w-z|^{-\alpha}+1]\, dw dz ds,\\
\label{eq:uest2}
Q_{1,1,1}^{x,y,t,t'} &\leq&  (8 + 3 K  2^{N_{\xi}\xi})^{2 \gamma}
\epsilon^{-2\delta_1\xi \gamma}
 \int_{0}^{t-\epsilon^{2}} (t-s)^{\xi \gamma}
\int_{A^{1}_{1}(x)}\int_{A^{1}_{1}(x)}
e^{\lambda \gamma |w|} e^{\lambda \gamma |z|} \\
\nonumber
&& \cdot \left| p_{t-s}(x-w)- p_{t'-s}(y-w)\right|\left| p_{t-s}(x-z)- p_{t'-s}(y-z)\right|
 [|w-z|^{-\alpha}+1]\, dw dz ds.
\end{eqnarray}
Note that the above integrals only become larger if we integrate over the
domain $[0,t]\times \IR^{2d},$ which we will do in the following.
We will use a version of the
factorization method first introduced
in \cite{DKZ87} to estimate them. Noting that
for $s\leq t$ and $0< a <1,$
\begin{equation}
\label{fac1}
\int_{s}^{t} (t-r)^{a-1} (r-s)^{-a} dr = \frac{\pi}{\sin(\pi a)},
\end{equation}
and that for $s \leq r \leq t,$
\begin{equation}
\label{fac2}
|p_{t-s}(x-w)-p_{t'-s}(y-w)|
\leq \int_{\IR^d} p_{r-s}(w'-w) \cdot |p_{t-r}(x-w') - p_{t'-r}(y-w')| dw'
\end{equation}
we obtain with (\ref{eq:uest1}),
\begin{eqnarray}
\label{fac3}
Q_{2,1,1}^{x,y,t,t'}
&\leq& C(a) \epsilon^{2(1-\delta_1)\xi \gamma}
\int_{0}^{t} \int_{0}^{t}
(t-r)^{a-1}  (t-r')^{a-1} \int_{\IR^{d}} \int_{\IR^{d}} J_{r,r'}^2(w',z') \\
\nonumber
& & \phantom{A}
|p_{t-r}(x-w')-p_{t'-r}(y-w')| \cdot  |p_{t-r'}(x-z')-p_{t'-r'}(y-z')| dw' dz' dr dr',
\end{eqnarray}
where
\begin{eqnarray}
\label{fac4}
J_{r,r'}^2(w',z')&:=& \int_{0}^{r \wedge r'} \int_{\IR^{d}}\int_{\IR^{d}}   
(r-s)^{-a} (r'-s)^{-a}
p_{r-s}(w'-w) p_{r'-s}(z'-z)\\
\nonumber
&&\phantom{AAAAAAAAAAAAAAAAAAAAAA}
\cdot [|w-z|^{-\alpha}+1]\,
dw dz ds,
\end{eqnarray}
where 
$J_{r,r'}^2(w',z') \leq J_{r,r'}^2(0,0)$  according to (\ref{eq:correst})
of Lemma \ref{lemma:correst}.  So we get
\begin{eqnarray}
\label{fac5}
Q_{2,1,1}^{x,y,t,t'}
&\leq& C(a)  \epsilon^{2(1-\delta_1)\xi \gamma} \int_{0}^{t}
\int_{0}^{t} (t-r)^{a-1}  (t-r')^{a-1} J_{r,r'}^2(0,0)   \\
\nonumber
& &
( \int_{\IR^{d}} |p_{t-r}(x-w')-p_{t'-r}(y-w')| dw' )
( \int_{\IR^{d}} |p_{t-r'}(x-z')-p_{t'-r'}(y-z')| dz' ) dr dr'.
\end{eqnarray}
The integrals in brackets can now be estimated with the help
of (\ref{eq:pdiffest}) in Lemma \ref{lemma:pdiffest}.  
Recall that $(t,x)\in Z_{K,N,\xi}$
and $|x-y|\le 2^{-N}$, so that $|x|\le K$, $|y|\le K+1$, and $t\le K$, and so
(\ref{eq:pdiffest}) implies
\begin{eqnarray}
\label{Qx211a}
\nonumber
Q_{2,1,1}^{x,y,t,t}&\leq& C(a,\alpha,d) \epsilon^{2(1-\delta_1)\xi \gamma} 
  |x-y|^{2\beta}
Q(t,a,0,\beta/2,\alpha)\\
&\le&C(\alpha,d,\beta,K)\epsilon^{2(1-\delta_1)\xi \gamma}
|x-y|^{2\beta}t^{1-\alpha/2-\beta}[1+t^{\alpha/2}].
\end{eqnarray}
Here we use $\beta < 1-\frac{\alpha}{2}$ and choose
$a\in (\beta/2,1-(\alpha/2))$ so that Lemma~\ref{lemma:Jest} may be applied in the last
line.  As
$t\le K$, (\ref{Qx211}) follows. Likewise we get for the time differences,
$\beta < 1-\frac{\alpha}{2},$ and ${\beta\over 2}<a<1-\alpha/2$, (use
Lemma~\ref{lemma:pdiffest} with $\beta/2$ in place of $\beta$),
\begin{eqnarray}
\nonumber
Q_{2,1,1}^{x,x,t,t'}&\leq& C(a,d) \epsilon^{2(1-\delta_1)\xi \gamma} 
 |t'-t|^{\beta}
Q(t,a,0,\beta/ 2,\alpha)\\
\label{Qt211a}
&\le& C(\beta,\alpha,d,K)
\epsilon^{2(1-\delta_1)\xi \gamma} |t'-t|^{\beta},
\end{eqnarray}
which is (\ref{Qt211}).

With an analogous calculation as in (\ref{fac1}) to (\ref{fac5}) except
now using (\ref{eq:uest2}) instead of (\ref{eq:uest1}),
we obtain that
\begin{eqnarray}
\label{fac6}
Q_{1,1,1}^{x,y,t,t'}
&\le& C(a) (8 + 3 K 2^{  N_{\xi} \xi})^{2\gamma} \epsilon^{-2\delta_1\xi \gamma }
\int_{0}^{t} \int_{0}^{t}
(t-r)^{a-1}  (t-r')^{a-1} \int_{\IR^{d}} \int_{\IR^{d}} J_{r,r'}^1(w',z') \\
\nonumber
& & \phantom{A}
|p_{t-r}(x-w')-p_{t'-r}(y-w')| \cdot  |p_{t-r'}(x-z')-p_{t'-r'}(y-z')| dw' dz' dr dr',
\end{eqnarray}
where
\begin{eqnarray}
\label{fac7}
J^{1}_{r,r'}(w',z') &=&
\int_{0}^{r \wedge r'} (t-s)^{\xi \gamma}
\int_{\IR^{d}}\int_{\IR^{d}}   (r-s)^{-a} (r'-s)^{-a}\\
\nonumber
& & \phantom{AAAAAAAAAA}p_{r-s}(w'-w) p_{r'-s}(z'-z)
\cdot e^{\lambda \gamma (|w| + |z|)} [|w-z|^{-\alpha}+1]\,
dw dz ds\\
\nonumber
&\leq & c_d^2 e^{\lambda^2(r+r')} e^{\lambda \gamma(|w'| + |z'|)}\int_{0}^{r \wedge r'}
(t-s)^{\xi
\gamma}
\int_{\IR^{d}}\int_{\IR^{d}}   (r-s)^{-a} (r'-s)^{-a}\\
\nonumber
& & \phantom{AAAAAAAAAAAAAA}p_{2(r-s)}(w'-w) p_{2(r'-s)}(z'-z) [|w-z|^{-\alpha}+1]\,
dw dz ds\\
\nonumber
&=: &   c_d^2 e^{\lambda^2(r+r')}e^{\lambda \gamma(|w'| + |z'|)}\tilde{J}_{r,r'}^{1}(w',z'),
\end{eqnarray}
where in the second inequality we have bounded $|w|\leq |w'| + |w'-w|$
(and likewise for $z$) and then used (\ref{eq:expest}).
Again $\tilde{J}_{r,r'}^{1}(w',z')\leq \tilde{J}_{r,r'}^{1}(0,0)$ independent of
$w'$ and $z'$ due to (\ref{eq:correst}) of Lemma \ref{lemma:correst}.
Hence, we get
\begin{eqnarray}
\label{fac8}
Q_{1,1,1}^{x,y,t,t'}
&\leq& C(d,a)e^{2\lambda^2K} (8 + 3 K 2^{  N_{\xi} \xi})^{2\gamma} \epsilon^{-2\delta_1\xi
\gamma }
  \int_{0}^{t}
\int_{0}^{t} (t-r)^{a-1}  (t-r')^{a-1} \tilde{J}_{r,r'}^1(0,0)   \\
\nonumber
& & \phantom{AAAAAAAAAA}
\times\left( \int_{\IR^{d}} |p_{t-r}(x-w')-p_{t'-r}(y-w')| e^{\lambda \gamma|w'|}dw' \right)
\\
\nonumber
& & \phantom{AAAAAAAAAA}
\times\left( \int_{\IR^{d}} |p_{t-r'}(x-z')-p_{t'-r'}(y-z')|
e^{\lambda \gamma|z'|} dz' \right) dr dr'.
\end{eqnarray}
And so, after a change of variables, the spatial differences are bounded by
\begin{eqnarray}
\nonumber
Q_{1,1,1}^{x,y,t,t}
&\leq& C(d,a)e^{2\lambda^2K} (8 + 3 K 2^{  N_{\xi} \xi})^{2\gamma}
e^{4\lambda^2K+2\lambda(K+1)+4\lambda}
\epsilon^{-2\delta_1\xi \gamma }
|x-y|^{2\beta'} 
Q(2t, a, \xi \gamma,\frac{\beta'}{2},\alpha)\\
\label{Qx111a}
&\leq& C(d,\xi\gamma,\beta',\alpha,K) (8 + 3 K 2^{  N_{\xi} \xi})^{2\gamma} 
\epsilon^{-2\delta_1\xi \gamma }
 |x-y|^{2\beta'}
\end{eqnarray}
if $\beta' < \xi \gamma + 1-\frac{\alpha}{2}$, $\beta' \le 1$, and $a$ is chosen in
$(\beta'/2,1-\alpha/2)\neq\emptyset$ (recall $\alpha<1$), according to (\ref{eq:pdiffest})
of Lemma
\ref{lemma:pdiffest} and Lemma
\ref{lemma:Jest} combined with (\ref{eq:expest}). This proves (\ref{Qx111}).
Similarly, for the time differences we obtain
\begin{eqnarray}
\nonumber
&&Q_{1,1,1}^{x,x,t,t'}\\
\nonumber
&&\leq
C(d,a) e^{2\lambda^2K}(8 + 3 K 2^{  N_{\xi} \xi})^{2\gamma}  e^{
4\lambda^2(K+1)+2\lambda(K+1)+4\lambda}  \epsilon^{-2\delta_1\xi \gamma } |t'-t|^{\beta'}
Q(0,0,0,2t, a,\xi \gamma,\frac{\beta'}{2},\alpha)\\
\label{Qt111a}
&&\leq C(d,\xi\gamma,\beta',\alpha,K) (8 + 3 K 2^{  N_{\xi} \xi})^{2\gamma}  
\epsilon^{-2\delta_1\xi \gamma } |t'-t|^{\beta'}
\end{eqnarray}
if again $\beta'< \xi \gamma + 1-\frac{\alpha}{2}$, $\beta'\le 1$ and $a$ is chosen as
above.  This shows (\ref{Qt111}).

Finally, we address the remaining case using (\ref{tubnd}) of Lemma
\ref{lemma:tubnd} to bound $\tilde{u}$ and Lemma \ref{lemma:correst}:
\begin{eqnarray}
\nonumber
Q^{x,t,t'}_{1,1} &\leq& C \epsilon^{2\gamma \xi(1-\delta_{1})}
\int_{t}^{t'} \int_{\IR^d} \int_{\IR^d} p_{t'-s}(x-w) p_{t'-s}(x-z)
[|w-z|^{-\alpha}+1]\, dw dz ds\\
\nonumber
&\leq&  C(\alpha,d) \epsilon^{2\gamma \xi(1-\delta_{1})} \int_{t}^{t'}
[(t'-s)^{-\frac{\alpha}{2}}+1] ds \\
\label{Qt11a}
&\leq& C(\alpha,d) \epsilon^{2\gamma \xi(1-\delta_{1})}
[|t'-t|^{1-\frac{\alpha}{2}}+|t'-t|],
\end{eqnarray}
and (\ref{Qt11}) follows as $|t'-t|\le 1$.
\end{proof}

\bigskip
Next, we consider all the terms for which $j=k=2.$ Here, we will use that
for $t \leq T_{K}$ we can bound
$|\tilde{u}(t,x)| \leq K e^{\lambda |x|}.$
\begin{lemma}
\label{lemma:Q22}
For $0<\beta< 1-\frac{\alpha}{2}$ we obtain for $i=1,2,$ and on $\{\omega:(t,x)\in
Z_{K,N,\xi}\}$,
\begin{eqnarray}
\label{Qxi22}
Q^{x,y,t,t}_{i,2,2}
&\leq& c_{\ref{Qxi22}}(d,\alpha,K)
 \exp(-\frac{1}{4} \epsilon^{-2\delta_{1}}(1-\beta)) |x-y|^{2\beta},\\
\label{Qti22}
Q^{x,x,t,t'}_{i,2,2}
&\leq& c_{\ref{Qti22}} (d,\alpha,K)
 \exp(-\frac{1}{4} \epsilon^{-2\delta_{1}}(1-{\beta\over 2})) |t'-t|^{\beta},\\
\label{Qt22}
Q^{x,t,t'}_{2,2}
&\leq& c_{\ref{Qt22}} (d,\alpha,\beta,K) 
\exp(-{1\over 2}\epsilon^{-2\delta_{1}}(1-\beta)) |t'-t|^{1-\frac{\alpha}{2}}.
\end{eqnarray}
\end{lemma}

\begin{proof} Recall $d((t,x),(t',y))\le \epsilon$.
For $i=1$ we are interested in  the case
$s \in [0,t-\epsilon^2]$ and $|x-w| > 2 \sqrt{t-s} \epsilon^{-\delta_{1}}.$
Since $|x-y|<\epsilon$ this implies that
$|y-w| \geq | \,|x-w| -|x-y|\,|> 2 \sqrt{t-s}\epsilon^{-\delta_{1}}-\epsilon
>\sqrt{t-s}\epsilon^{-\delta_{1}}.$ Furthermore, $t'-s=t'-t+t-s\leq \epsilon^2 +t-s
\leq 2(t-s). $ This implies
\begin{equation}
\label{expbound}
\exp(-\frac{|x-w|^{2}}{4(t'-s)}) \vee \exp(-\frac{|y-w|^{2}}{4(t'-s)})
\leq \exp(-\frac{|x-w|^{2}}{8(t-s)}) \vee \exp(-\frac{|y-w|^{2}}{8(t-s)})
\leq \exp(-\frac{1}{8} \epsilon^{-2\delta_{1}}).
\end{equation}
Therefore, for $v=x$ or $v=y$ and $r=t$ or $r=t',$
\begin{eqnarray}
\label{eq:pexpbound}
p_{r-s}(v-w)\leq 2^{\frac{d}{2}}  \exp(-\frac{1}{8} \epsilon^{-2\delta_{1}}) p_{2(r-s)}(v-w).
\end{eqnarray}
Using this we obtain for any $\beta \in (0,1)$  by
applying H\"older's inequality that $Q^{x,y,t,t'}_{1,2,2}$ is bounded by
\begin{eqnarray}
\nonumber
&&\int_{0}^{t-\epsilon^2}
\Big(\int_{A^{1,s}_{2}(x)} \int_{A^{1,s}_{2}(x)} \left( p_{t-s}(x-w)+ p_{t'-s}(y-w)\right)
\left( p_{t-s}(x-z)+ p_{t'-s}(y-z)\right) \\
\nonumber
&&\phantom{\cdot \Big(\int_{A^{1,s}_{2}}
\int_{A^{1,s}_{2}} p_{t-s}(x-w)  AAAAAAA}
\cdot [|w-z|^{-\alpha}+1]  |\tilde u(s,w)|^{\frac{\gamma}{1-\beta}}|\tilde u(s,z)|^{\frac{\gamma}{1-\beta}}
dw dz\Big)^{1-\beta}\\
\nonumber
&&\phantom{}
\cdot \Big(\int_{A^{1,s}_{2}(x)} \int_{A^{1,s}_{2}(x)} \left| p_{t-s}(x-w)- p_{t'-s}(y-w)\right|
\left| p_{t-s}(x-z)- p_{t'-s}(y-z)\right|\\
\nonumber
&&\phantom{\cdot \Big(\int_{A^{1,s}_{2} }
\int_{A^{1,s}_{2}} \left| p_{t-s}(x-w)- p_{t'-s}(y-w)\right|  AAAAAAAAA}\cdot [|w-z|^{-\alpha}+1] dw
dz\Big)^{\beta} ds\\
&\leq& C(d,\alpha) e^{2\lambda^2K}K^{2 \gamma}e^{2\lambda \gamma (K+1)}
\exp(-\frac{1}{4} \epsilon^{-2\delta_{1}}(1-\beta))
\int_{0}^{t-\epsilon^2} [(t-s)^{-\frac{\alpha }{2}(1-\beta)}+1]\\
\nonumber
&&\phantom{}
\Big(\int_{\IR^d} \int_{\IR^d} \left| p_{t-s}(x-w)- p_{t'-s}(y-w)\right|
\left| p_{t-s}(x-z)- p_{t'-s}(y-z)\right| [|w-z|^{-\alpha}+1] dw dz\Big)^{\beta} ds.
\end{eqnarray}
Here, we have used that $\tilde{u}(t,x) \leq K e^{\lambda |x|}$ and (\ref{eq:pexpbound})
as well as  (\ref{eq:correstlambda}) of Lemma \ref{lemma:correst}. We have also used the fact that
$e^{\lambda \gamma |x|}$ and $e^{\lambda \gamma |y|}$ are both bounded by
$e^{\lambda \gamma (K+1)}$ since $|x|<K$ and $|x-y|<\epsilon<1.$

\noindent
Using (\ref{eq:spacecorrest}) of Lemma \ref{lemma:pdiffest} to estimate the integral
in parentheses when $t=t'$ we obtain
\begin{eqnarray}
\nonumber
Q^{x,y,t,t}_{1,2,2}
&\leq& C(d,\alpha)e^{2\lambda^2K} (K e^{\lambda (K+1)})^{2}
\exp(-\frac{1}{4} \epsilon^{-2\delta_{1}}(1-\beta)) |x-y|^{2\beta}\\
\nonumber
&&\phantom{C(d,\alpha)e^{2\lambda^2K} (K e^{\lambda (K+1)})^{2}}\cdot
\int_{0}^{t-\epsilon^2} (t-s)^{-\frac{\alpha }{2}-\beta}+(t-s)^{-\beta} ds\\
\label{Qx122a}
&\leq& C(d,\alpha)K^3 e^{2(\lambda^2+\lambda)(K+1)} 
 \exp(-\frac{1}{4} \epsilon^{-2\delta_{1}}(1-\beta)) |x-y|^{2\beta},
\end{eqnarray}
provided that $\beta< 1-\frac{\alpha}{2}$, showing (\ref{Qxi22}) for i=1.
Likewise, using (\ref{eq:timecorrest}) of Lemma \ref{lemma:pdiffest}
for the time differences when $x=y$ implies (\ref{Qti22}) for i=1
if again $0<\beta< 1-\frac{\alpha}{2}$--here we replace $\beta$ with $\beta/2$ in the above.

For $i=2$ and $Q_{2,2}^{x,t,t'}$ we will proceed analogously. We merely have to establish
(\ref{expbound}) in the case:
$s \in [t-\epsilon^2,t']$ and $|x-w| > 2 \epsilon^{1-\delta_{1}}.$
Since $|x-y|<\epsilon$ this implies now
$|y-w| \geq | |x-w| -|x-y||> 2 \epsilon^{1-\delta_{1}}-\epsilon
>\epsilon^{1-\delta_{1}}.$ Furthermore, $t'-s=t'-t+t-s\leq \epsilon^2 +t-s
\leq 2 \epsilon^{2}. $ From this the bound (\ref{expbound}) follows and we obtain
immediately (\ref{Qxi22}) and (\ref{Qti22}) for $i=2$ provided that
$\beta< 1-\frac{\alpha}{2}.$

Lastly, we obtain with the
help of (\ref{expbound}) (verified above) and (\ref{eq:correstlambda}) of Lemma
\ref{lemma:correst},
\begin{eqnarray}
\nonumber
Q^{x,t,t'}_{2,2} &\leq& C(d,\beta) K^{2\gamma} \exp(-\frac{1}{2} \epsilon^{-2\delta_{1}}(1-\beta))\\
\nonumber
& &\phantom{Q^{x,t,t'}_{2,2} }
\cdot   \int_{t}^{t'} \int_{\IR^{d}} 
\int_{\IR^{d}} e^{\lambda
\gamma (|w|+|z|)} p_{\frac{t'-s}{\beta}}(x-w)  p_{\frac{t'-s}{\beta}}(x-z)
[|w-z|^{-\alpha}+1]\, dw dz ds\\
\nonumber
&\leq& C(d,\beta, \alpha)  K^{2\gamma}e^{2\lambda^2(K+1)/\beta}e^{2\lambda(K+1)}
\exp(-\frac{1}{2} \epsilon^{-2\delta_{1}}(1-\beta))
\int_{t}^{t'} [(\frac{2}{\beta}(t'-s))^{-\frac{\alpha}{2}}+1]\, ds \\
\label{Qt22a}
&\leq& C(d,\alpha, \beta,K)) \exp(-\frac{1}{2}
\epsilon^{-2\delta_{1}}(1-\beta)) |t'-t|^{1-\frac{\alpha}{2}},
\end{eqnarray}
which is (\ref{Qt22}) and hence completes the proof.
\end{proof}

It remains to consider the ``mixed terms'' for which $j=2$ and $k=1$ or vice versa.
Say $j=2.$ In this case (\ref{expbound}) holds for the exponential in the $w$ integral,
and we can bound the exponential in the $z$ integral by one. Otherwise we follow the same
steps as in Lemma \ref{lemma:Q22} treating the case $j=k=2.$
In this manner, we obtain the same bounds as in (\ref{Qxi22})
to (\ref{Qt22}) with the only difference that
$\exp(-\frac{1}{4} \epsilon^{-2\delta_{1}}(1-\beta))$
is replaced by $\exp(-\frac{1}{8} \epsilon^{-2\delta_{1}}(1-\beta))$ and
$\exp(-\frac{1}{2} \epsilon^{-2\delta_{1}}(1-\beta))$
by $\exp(-\frac{1}{4} \epsilon^{-2\delta_{1}}(1-\beta)).$

We are now taking the estimates (\ref{Qx211}), (\ref{Qx111}) and (\ref{Qxi22})
together with those for the mixed terms and choose $\beta=1-\frac{\alpha}{2}-\delta'$,
respectively, $\beta'=1-\frac{\alpha}{2}-\delta'+\xi \gamma <1$ (by (\ref{chdelta'})) 
in those
estimates. This shows that for $(t,x)\in Z_{K,N,\xi},\ |x-y|<\epsilon=2^{-N}$
and $N\ge N_1$,
\begin{eqnarray}
\nonumber
Q^{x,y,t,t} &\leq& C(K) |x-y|^{2(1-\frac{\alpha}{2}-\delta')} \Big[
 \epsilon^{2 (1- \delta_{1}) \xi \gamma }
+ (8+ 3K 2^{ N_{\xi} \xi})^{2\gamma}  
\epsilon^{-2\delta_1\xi \gamma } |x-y|^{2\xi \gamma}\\
\nonumber
& &\phantom{C(K) |x-y|^{2(1-\frac{\alpha}{2}-\delta')}}
+ \exp( -\frac{1}{4} \epsilon^{-2\delta_{1}} (\frac{\alpha}{2}+\delta' )) +
\exp( -\frac{1}{8} \epsilon^{-2\delta_{1}} (\frac{\alpha}{2}+\delta')) \Big] \\
\label{Qx}
&\leq& C(K) |x-y|^{2(1-\frac{\alpha}{2}-\delta')} \Big[ \epsilon^{2 (1- \delta_{1})\xi
\gamma} 2^{2  N_{\xi} \xi \gamma} + \exp(-\frac{\alpha}{16}  \epsilon^{-2\delta_{1}}) \Big].
\end{eqnarray}

\noindent
We have the analogous bounds
for $Q^{x,x,t,t'}+Q^{x,t,t'}$ with the help of
(\ref{Qt211}), (\ref{Qt111}), (\ref{Qt11}), (\ref{Qti22}), and (\ref{Qt22});
just replace $|x-y|^{2}$ with $|t'-t|$ and use $ |t'-t| < \epsilon^{2}$.
We deduce that for $N\ge N_1$ and $(t,x)\in Z_{K,N,\xi}$,
\begin{equation}
\label{eq:Q_2bound}
Q^{x,x,t,t'}+Q^{x,t,t'}\leq C(K) |t'-t|^{1-\frac{\alpha}{2}-\delta'}
\Big[ \epsilon^{2 (1- \delta_{1})\xi \gamma}
2^{2  N_{\xi} \xi \gamma} + \exp(-\frac{\alpha}{16}  \epsilon^{-2\delta_{1}}) \Big].
\end{equation}
\noindent
We can finally conclude that in (\ref{eq:H3}), $P_{3}=P_{4}=0$ if
\begin{equation}
\label{eq:Hoelderimp}
C(K) \Big[ \epsilon^{2 (1- \delta_{1})\xi \gamma}
2^{2  N_{\xi} \xi \gamma} + \exp(-\frac{\alpha}{16}  \epsilon^{-2\delta_{1}}) \Big]
< \epsilon^{2p}.
\end{equation}
For this it is sufficient that
\begin{eqnarray}
\label{eq:small1}
C(K) \epsilon^{2 (1- \delta_{1})\xi \gamma}
2^{2  N_{\xi} \xi \gamma} &<&\frac{1}{2} \epsilon^{2p},\\
\label{eq:small2}
C(K) \exp(-\frac{\alpha}{16}  \epsilon^{-2\delta_{1}}) &<&\frac{1}{2} \epsilon^{2p}.
\end{eqnarray}
Since (\ref{eq:small1}) is equivalent to
$2C(K) <  2^{ 2N [(1- \delta_{1})\xi \gamma-p] - 2N_{\xi} \xi \gamma}$,
it suffices to choose $\delta_1>0$ small enough so that
$(1- \delta_{1})\xi \gamma-p >0$
(which is possible since $\xi \gamma >p$) and then to assume $N\geq
[C_0(\xi,\delta_1) N_{\xi}]\in\IN$ as well as $N\geq N_{0}(K,\xi,\delta_{1},p)
\in \IN$ deterministic so that both (\ref{eq:small1}) and (\ref{eq:small2}) hold.
Note that the constants depend ultimately on $\xi, \xi_{1}$ and $K.$
Hence, (\ref{eq:H3}), (\ref{P1est}) and (\ref{P2est}) imply that if
$N_2(\omega,\xi,\xi_1,K)=\Bigl[{5N_\xi\over
\delta_1}\Bigr]\vee [C_0(\xi,\delta_1) N_{\xi}] \vee  N_{0}(K,\xi,\delta_{1},p)$, then for
$d((t,x),(t',y))\le 2^{-N}$, $t\le t'$,
\begin{eqnarray}
\nonumber
& &\IP\Big(|\tilde u(t,x)-\tilde u(t,y)|\ge |x-y|^{1-\frac{\alpha}{2}-\delta}2^{-Np},\ (t,x)\in
Z_{K,N,\xi}, N\ge N_2 \Big )\\
\nonumber
& &\qquad+\IP \Big(|\tilde u(t',x)-\tilde u(t,x)|\ge |t'-t|^{\frac{1}{2}
(1-\frac{\alpha}{2}-\delta)}2^{-Np},
\ (t,x)\in Z_{K,N,\xi}, t'\le T_K, N\ge N_2 \Big)\\
\label{keybnd}
& &\le
c_{\ref{P1est}}(\exp\Bigl(-c'_{\ref{P1est}}|x-y|^{-\delta''}\Bigr)
+\exp\Bigl(-c'_{\ref{P1est}}|t'-t|^{-\frac{\delta''}{2}}\Bigr)).
\end{eqnarray}
\noindent
Now let $e_{l}$ be the $l^{\hbox{th}}$ unit vector in $\IR^d$ and set
\begin {eqnarray*}
M_{n,N,K}&=& \max\{\sum_{l=1}^{d}|\tilde u(j2^{-2n},(z+e_{l})2^{-n})-\tilde u(j2^{-2n},z2^{-n})|\\
& &\qquad +|\tilde u((j+e)2^{-2n},z2^{-n})-\tilde u(j2^{-2n},z2^{-n})|:\\
& &\qquad |z|\le K2^n, (j+e)2^{-2n}\le T_K, j\in \IZ_+, z\in\IZ^{d}, \\
& &\qquad e\in\{1,2,3\},  (j2^{-2n},z2^{-n})\in Z_{K,N,\xi}\}.
\end{eqnarray*}
(\ref{keybnd}) implies that if
\[A_N=\{\omega:\hbox{ for some }n\ge N,\ M_{n,N,K}\ge (d+1)\cdot
2^{-n(1-\frac{\alpha}{2}-\delta)}2^{-Np},\ N\ge N_2\},
\]
then for some fixed constants $C(d),c_1,c_2>0$, 
\begin{eqnarray*}
\IP(\cup_{N'\ge N}A_{N'})&\le& C(d)\sum_{N'=N}^\infty\sum_{n=N'}^\infty
K^{d+1} 2^{(d+2)n} e^{-c_1 2^{n\delta''}}\\
&\le&  C(d)K^{d+1} \eta_N,\\
\end{eqnarray*}
where $\eta_N =e^{-c_{2} 2^{N \delta''}}.$
Therefore $N_3(\omega)=\min\{N\in\IN:\omega\in A_{N'}^c\ \hbox{for all }N'\ge
N\}<\infty$ a.s. and in fact
\begin{eqnarray}
\label{N3bnd}
\IP(N_3> N)=\IP(\cup_{N'\ge N}A_{N'})\le  C(d)K^{d+1}\eta_N.
\end{eqnarray}
Choose $m \in \IN$ with $m > \log_{2}( 3+ \sqrt{d})$ and assume
$N\ge (N_3+m)\vee N_2$.  Let $(t,x)\in Z_{K,N,\xi}$, $d((t',y),(t,x))\le
2^{-N}$, and $t'\le T_K$.  For $n\ge N$ let $t_n\in 4^{-n}\IZ_+$ and $x_{n,i}\in
2^{-n}\IZ$ ($i=1,\dots,d$) be the unique points so that $t_n\le t<t_n+4^{-n},$
$x_{n,i}\le x_i<x_{n,i}+2^{-n}$
for $x_{i}\ge 0$ and $x_{n,i}-2^{-n}<x_i\le x_{n,i}$ if $x_i<0$.
Similarly define $t'_n$ and
$y_n$ with $(t',y)$ in place of $(t,x)$.  Choose $(\hat{t},\hat{x})$ as in the definition of
$Z_{K,N,\xi}$ (recall $(t,x)\in Z_{K,N,\xi}$).  If $n\ge N$, then
\begin{eqnarray*}
d((t'_n,y_n),(\hat{t},\hat{x}))&\le&
d((t'_n,y_n),(t',y))+d((t',y),(t,x))+d((t,x),(\hat{t},\hat{x}))\\
 &\le& \sqrt{|t'_n-t'|}+|y-y_n|+2^{-N}+2^{-N}\\
&<& (3+ \sqrt{d})2^{-N}< 2^{m-N}.
\end{eqnarray*}
Therefore $(t'_n,y_n)\in Z_{K,N-m,\xi}$, and similarly (and slightly more simply)
$(t_n,x_n)\in Z_{K,N-m,\xi}$.  Our definitions imply
that $t_N$ and $t'_N$ are equal or adjacent in $4^{-N}\IZ_+$ and similarly for
the components of $x_N$
and $y_N$ in $2^{-N}\IZ_+.$  This, together with the continuity of $\tilde u$, the 
triangle inequality,
and our lower bound on $N$ (which shows $N-m\ge N_3$), implies
\begin{eqnarray*}
|\tilde u(t,x)-\tilde u(t',y)|&\le& |\tilde u(t_N,x_N)-\tilde
u(t'_N,y_N)|\\
& &\qquad+\sum_{n=N}^\infty|\tilde u(t_{n+1},x_{n+1})-\tilde u(t_n,x_n)|
+|\tilde u(t'_{n+1},y_{n+1})-\tilde u(t'_n,y_n)|\\
&\le& M_{N,N-m,K}+\sum_{n=N}^\infty 2M_{n+1,N-m,K}\\
&\le& 4\sum_{n=N}^\infty (d+1)\cdot2^{-n(1-\frac{\alpha}{2}-\delta)}2^{-(N-m)p}\\
&\le& c_0(d,p)2^{-N(1-\frac{\alpha}{2}       -\delta+p)}\\
&\le& 2^{-N\xi_1}.
\end{eqnarray*}
The last line is valid for $N\ge N_4$ because $1-\frac{\alpha}{2}-\delta+p>\xi_1$ by
(\ref{chp}). Here
$N_4$ is deterministic and may depend on $p,\xi_1,\delta,c_0$ and hence ultimately on
$\xi,\xi_1$.  This proves the required result with
\[N_{\xi_1}(\omega)=\max(N_3(\omega)+m,\Bigl[\frac{5N_\xi(\omega)}
{\delta_1}\Bigl],[C_{0}(\xi,\delta_{1}) N_{\xi}],N_{0}\vee N_4).
\]
Therefore, if $R=5/\delta_1 \vee C_{0}(\xi,\delta_{1})$
and $N\geq N(K):=N_{0}\vee N_4$ (deterministic), (\ref{N3bnd})
implies that
\[\IP(N_{\xi_1}\ge N)\le\IP(N_3\ge N-m)+2 \IP(N_\xi\ge N/R)\le
c(d) K^{d+1} \eta_{N-m}+2 \IP(N_{\xi}\ge N/R),\]
which gives the required probability bound (\ref{Nbnd}).

\bigskip

\section{Appendix (Proof of Theorems~\ref{thm:nonLipschitzexistence},~\ref{prop})}

In this appendix, we briefly describe the construction of solutions to (\ref{heatecolint})
with colored noise and non-Lipschitz coefficients. We start by citing the
following result which states necessary conditions for the existence of solutions to
(\ref{heatecolint}) with Lipschitz coefficients and bounded initial conditions
(see Dalang \cite{rD99}):
\begin{thm}
\label{thm:Dalang}
Let $u_{0}$ be measurable and bounded and let $\sigma$ be a Lipschitz continuous function.
Assume that $(A)_{\eta}$ holds for $\eta=1.$ Then there exists a pathwise unique
solution $u$ to (\ref{heatecolint}) which is also a strong solution.
The process $u$ satisfies a uniform moment bound:
For any $T>0,$ and $p \in [1,\infty),$
\begin{equation}
\label{eq:uniformmoments}
\sup_{0 \leq t \leq T} \sup_{x \in \IR^{d}}
\IE\left( |u(t,x)|^{p} \right) < \infty.
\end{equation}
\end{thm}
\noindent
We would like to remark that the original theorem of Dalang \cite{rD99}
stipulates that the noise be spatially homogeneous. However, it is not hard to
see that all that is needed is that it be bounded by an appropriate spatially
homogeneous term in the sense of condition $(A)_{\eta}.$

Denote $L^{\infty}_{tem}=\left\{u: \;\mbox{ess\,sup}_{x \in \IR^{d}}
 |u(x)| e^{-\lambda |x|} <\infty\hbox{ for all }\lambda>0\right\}$.  Here the
$\mbox {ess\,sup}$ is of course with respect to Lebesgue measure.

We introduce some frequently-used notation. For any
function 
$v:\IR_{+} \times \IR^{d} \rightarrow \IR$ and stopping time $\tau$, we set 
\begin{eqnarray}
J^{a-1}v(t,x)&=&\frac{\sin(\pi a)}{\pi}
\int_{0}^{t}\int_{\IR^{d}} (t-s)^{a-1}
p_{t-s}(x-y)v(s,y) dy ds,
\end{eqnarray}
as well as
\begin{eqnarray}
J^{\tau}_{a} v(t,x) &=& \int_{0}^{t}\int_{\IR^{d}} 1(s\leq \tau)(t-s)^{-a}
p_{t-s}(x-y) \sigma(v(s,y)) W(dyds).
\end{eqnarray}
The stochastic Fubini Theorem implies
\begin{eqnarray}
\label{equt:111}
J^{a-1}J_{a}^{\tau} v(t,x)
= \int_{0}^{t} \int_{\IR^{d}}1(s\leq \tau) p_{t-s}(x-y)
\sigma(v(s,y)) W(dyds).
\end{eqnarray}
We will use the notation $J_{a}v(t,x)=J_{a}^{t}v(t,x)$, when $\tau=t$ in the above. 
Also set 
\begin{equation*}
G^{\tau}_{\lambda,p}v(t,x):= \IE \left(|v(t,x)|^p 1(t\leq\tau)
e^{-\lambda |x|} \right),
\end{equation*}
and again $G_{\lambda,p}v(t,x)\equiv G^{t}_{\lambda,p}v(t,x)$ whenever $\tau=t$. 
\begin{lemma}
\label{lem:2} 
Let $\sigma$ be a
continuous function satisfying the growth condition
\begin{eqnarray}
\label{equt:103} 
  |\sigma(u)|\leq c_{\ref{equt:103}}(1+|u|).
\end{eqnarray}
 Assume that  $(A)_\eta$ holds for some $\eta \in [0,1)$ and let $a<(1-\eta)/2$.   
Let $v:\Omega\times\IR_+\times\IR^d\to\IR$ be ${\cal P}({\cal F}_\cdot)\times
{\cal B}(\IR^d)$-measurable (${\cal P}({\cal F}_\cdot)$ is the $({\cal
F}_t)$-predictable
$\sigma$-field). Then for any
$T,\lambda>0$, $p\geq 2$, and stopping time $\tau$,  
\begin{equation}
\label{equt:122}
\sup_{0 \leq s \leq t} \sup_{x \in \IR^{d}}\IE \left(
 |J_{a}^{\tau} v(t,x)|^{p} e^{-\lambda |x|} \right) \leq
C(T,\lambda,p)c_{\ref{equt:103}}^{p}
\sup_{0\leq s\leq t}
 \sup_{x\in\IR^d} 
\left(1+  G^{\tau}_{\lambda,p}v(s,x)
\right),\;\; \forall t\leq T. 
\end{equation}
\end{lemma}
\begin{proof}
First fix arbitrary $p\geq 2$ and $x\in \IR^d$. 
Then, 
using the growth condition on $\sigma$ as well as
Burkholder's inequality and $|k(x,y)|\leq c_{\ref{eq:colorcond}} \tilde{k}(x-y)$ we
get
\begin{eqnarray}
\nonumber
& &\IE \Big( |J_{a}^{\tau}v(t,x)|^{p}  \Big)\\
\nonumber
&\leq& C\, c_{\ref{equt:103}}^p
\IE \Big( \big(\int_{0}^{t}\int_{\IR^{d}}\int_{\IR^{d}}
(t-s)^{-2 a} p_{t-s}(x-y) p_{t-s}(x-z) \tilde{k}(y-z) \\
\nonumber
& &\phantom{AAAAAAAAAAAAAAAAAAA} \cdot (1+
 1(s\leq \tau)|v(s,y)|)(1+ 1(s\leq \tau)|v(s,z)|)
dy dz ds \big)^{\frac{p}{2}} \Big) \\
\nonumber
&\leq& C\,  c_{\ref{equt:103}}^p\left(\int_{0}^{t}
(t-s)^{-2 a} \int_{\IR^{d}} \int_{\IR^{d}}
p_{t-s}(y) p_{t-s}(z) \tilde{k}(y-z) dy dz ds \right)^{\frac{p}{2}-1}\\
\nonumber
& & \phantom{C} \cdot \Big(\int_{0}^{t}
(t-s)^{-2 a} \int_{\IR^{d}} \int_{\IR^{d}} p_{t-s}(y)  p_{t-s}(z)
\tilde{k}(y-z)\\
\nonumber
& & \phantom{AAAAAAA}
\cdot \IE \Big((1+ 1(s\leq \tau)|v(s,y-x)|)^{\frac{p}{2}}
(1+ 1(s\leq \tau)|v(s,z-x)|)^{\frac{p}{2}} \Big) dy dz  ds\Big).
\end{eqnarray}
Apply H\"older's inequality to the expected value in this
expression and shift variables to bound it by
\begin{eqnarray}
\nonumber
\lefteqn{\IE \Big( (1+ 1(s\leq \tau)|v(s,y-x)|)^{p} \Big)^{\frac{1}{2}}
\IE \Big((1+ 1(s\leq \tau)|v(s,z-x)|)^{p} \Big)^{\frac{1}{2}}}\\
\label{EHoelder}
&\leq&  C(\lambda,p) e^{\frac{\lambda}{2} (|y|+|z|)+\lambda|x|} (1+\sup_{\tilde{z}\in \IR^d} G^{\tau}_{\lambda,p}
v(s,\tilde{z})).
\end{eqnarray}
Hence, we arrive at
\begin{eqnarray}
\nonumber
& &\IE \Big(  |J_{a}^{\tau} v(t,x)|^{p}
e^{-\lambda |x|} \Big)\\
\nonumber
&\leq& C(\lambda,p) c_{\ref{equt:103}}^p \left(\int_{0}^{t}
(t-s)^{-2 a} \int_{\IR^{d}}\int_{\IR^{d}}
 p_{t-s}(y) p_{t-s}(z) \tilde{k}(y-z) dy dz ds \right)^{\frac{p}{2}-1}\\
\nonumber
& & \cdot \left(\int_{0}^{t}
(t-s)^{-2 a} (\int_{\IR^{d}} \int_{\IR^{d}} e^{\frac{\lambda}{2}(|y|+|z|)}
p_{t-s}(y)  p_{t-s}(z) \tilde{k}(y-z) dy dz)
(1+ \sup_{\tilde{z}\in \IR^d} G^{\tau}_{\lambda,p}v(s,\tilde{z})) ds\right)\\
&\leq&
\nonumber
C(\lambda,p) c_{\ref{equt:103}}^p \left( \int_{0}^{t} f(\frac{s}{2})
ds\right)^{\frac{p}{2}-1}
\left(\int_{0}^{t} f(t-s) (1+\sup_{\tilde{z}\in \IR^d} G^{\tau}_{\lambda,p}v(s,\tilde{z})) ds\right)\\
\label{equt:113}
&\leq& C(T,\lambda,p) c_{\ref{equt:103}}^p \sup_{0\leq s\leq
t}\sup_{\tilde{z}\in \IR^d} (1+G^{\tau}_{\lambda,p}v(s,\tilde{z})),\;\; 
 \forall t\leq T, \; x\in \IR^d, 
\end{eqnarray}
where
\begin{equation*}
f(r) = r^{-2 a}
(\int_{\IR^{d}} \int_{\IR^{d}}
p_{2r}(y)  p_{2r}(z) \tilde{k}(y-z) dy dz).
\end{equation*}
Here, we have used that
$e^{\frac{\lambda}{2}|y|}p_{t}(y)\leq C(T,\lambda) p_{2t}(y)$
for $t \leq T,$ see (\ref{eq:expest}). We have also used the fact
that $f$ is integrable on $[0,T]$ for $a< \frac{1-\eta}{2}$
(cf. proof of Lemma 2.2 of \cite{SS02}).
This 
proves~(\ref{equt:122}) for all $p\geq 2$.
\end{proof}

\begin{lemma}
\label{lem:1} 
Let $u_{0}\in L^{\infty}_{tem}$ and  let $\sigma$ be a
continuous function satisfying the growth condition~(\ref{equt:103}).
 Assume that  $(A)_\eta$ holds for some $\eta \in [0,1).$  
If  $u$  is any solution to~(\ref{heatecolint}) such that 
\begin{eqnarray}
\sup_{0\leq t\leq T}\sup_{x\in\IR^d}\IE \left( \left|u(t,x)
 \right|^p e^{-\lambda |x|}\right)&<&\infty 
,\;\; \forall T>0, p>0, \lambda>0, 
\end{eqnarray}
then for any $T,p,\lambda>0$, there exists $\tilde{p}\geq p$ such that 
\begin{equation}
\label{equt:121}
\IE \left(\sup_{0 \leq t \leq T} \sup_{x \in \IQ^{d}}
 |u(t,x)|^{p} e^{-\lambda |x|} \right) \leq C_{T,\lambda,p}
(c_{\ref{equt:103}},\left\|u_{0}
\right\|_{\frac{\lambda}{p},\infty})\left(1+ 
\sup_{0\leq t\leq T}\sup_{x\in\IR^d} G_{\frac{\lambda}{2},\tilde{p}}u(t,x)
\right),
\end{equation}
where $C_{T,\lambda,p}(\cdot,\cdot)$ is bounded on the compacts of
 $\IR_+\times\IR_+$. 
\end{lemma}
\begin{proof}
\begin{eqnarray}
\nonumber
\lefteqn{\IE \left( \sup_{0\leq t\leq T}\sup_{x\in \IQ^d} \left|u(t,x)
 \right|^p e^{-\lambda |x|}\right)
\leq
C\IE \left( \sup_{0\leq t\leq T}\sup_{x\in \IQ^d}
\left|\int_{\IR^d} p_t(x-y) u_0(y) dy\right|^p e^{-\lambda |x|} \right)}\\
\label{equt:118}
& &+ C\IE \left( \sup_{0\leq t\leq T}\sup_{x\in \IQ^d}
\left|\int_0^{t} \int_{\IR^d} p_{t-s}(x-y) 
 \sigma(u(s,y)) W(dyds)\right|^p
e^{-\lambda |x|}  \right).
\end{eqnarray}
The first term on the right hand side of~(\ref{equt:118}) is bounded by
\begin{eqnarray}
\nonumber
\lefteqn{\sup_{0\leq t\leq T}\sup_{x\in \IR^d}\left|\int_{\IR^{d}} p_t(x-y) |u_0(y)|^{p} dy e^{-\lambda |x|}\right|}
\\
\nonumber
&\leq& C(T,\lambda) \left\|u_{0}
\right|^p_{\frac{\lambda}{p},\infty} \sup_{0\leq t\leq T}\sup_{x\in
\IR^d}\left|
 \int_{\IR^{d}} p_t(x-y)  e^{\lambda |y|} dy e^{-\lambda |x|}\right|
\\
\label{equt:119}
&\leq&
 C(T,\lambda,p)\left\|u_{0}
 \right\|^p_{\frac{\lambda}{p},\infty}.
\end{eqnarray}
In this calculation we have used Jensen's Inequality as well as the fact that
\begin{equation}
\label{equt:104}
\int_{\IR^{d}}p_{t}(x-y) e^{\lambda |y|} dy \leq C(T,\lambda) e^{\lambda |x|}
\end{equation}
for $t \leq T$ and $ \lambda \in \IR$ (see Lemma 6.2 of \cite{tS94}).

We bound the second  term on the right hand side of~(\ref{equt:118})
  with the help of the factorization
method of \cite{DKZ87} (compare (\ref{fac1}) and (\ref{fac2})).
Let $0<a<(1-\eta)/2$ and choose arbitrary $p^{*}> \frac{1+\frac{d}{2}}{a}>2$. 
Assume that 
$p\geq p^{*}$. Recall $\Vert v\Vert_{\lambda,p}=\Bigl[\int
|v(x)|^pe^{-\lambda |x|}\,dx\Bigr]^{1/p}$.  Use~(\ref{equt:111}) and apply 
H\"older's inequality to get 
\begin{eqnarray}
\nonumber
\lefteqn{\IE\left( \sup_{t \leq T}
\sup_{x \in \IQ^{d}}\left| \int_{0}^{t} \int_{\IR^{d}}p_{t-s}(x-y)
\sigma(u(s,y)) W(dyds)\right|^{p}e^{-\lambda |x|}\right)}
\\
\nonumber
&=& \IE \left(  \sup_{t \leq T}
\sup_{x \in \IQ^{d}} |J^{a-1}J_{a}u(t,x) |^{p}
e^{-\lambda |x|} \right)\\
\nonumber
&\leq& C
\IE \left(  \sup_{t \leq T} \sup_{x \in \IQ^{d}}| \int_{0}^{t}  (t-s)^{a-1}
\left( \int_{\IR^{d}} p_{t-s}(x-y)
e^{\frac{\lambda}{2} |y|}
\cdot |J_{a}u(s,y)|^{\frac{p}{2}}e^{-\frac{\lambda}{2} |y|} dy
\right)^{\frac{2}{p}}  ds |^{p} e^{-\lambda |x|}  \right) \\
\nonumber
&\leq& C
\IE \left( \sup_{t \leq T} \sup_{x \in \IQ^{d}}|\int_{0}^{t}  (t-s)^{a-1}
\left( \int_{\IR^{d}} p_{t-s}(x-y)^{2} e^{\lambda |y|} dy \right)^{\frac{1}{p}}
\cdot || J_{a} u_{s}||_{\lambda,p}   ds
|^{p} e^{-\lambda |x|}  \right) \\
\nonumber
&\leq& C(T,\lambda) \IE \left( \sup_{t \leq T}
\left( \int_{0}^{t}  (t-s)^{a-1-\frac{d}{2p}}
\cdot || J_{a}u_{s}||_{\lambda,p} ds \right)^{p}   \right)\\
\label{equt:115}
&\leq& C(T, \lambda) \left(\int_{0}^{T}
s^{(a-1-\frac{d}{2p})\frac{p}{p-1}} ds \right)^{p-1}
\cdot  \int_0^T\IE\left(|| J_{a}u_{s}||_{\lambda,p}^{p} \right)\,ds.
\end{eqnarray}
Here, we have also used 
(\ref{equt:104})
and $p_{t}(x) \leq C t^{-\frac{d}{2}}$.
Lemma~\ref{lem:2}
implies
\begin{eqnarray}
\nonumber
\IE \Big( \int_{\IR^d} |J_{a} u(t,x)|^{p}
e^{-\lambda |x|}\,dx  \Big)
&\leq& C(T,\lambda,p)c_{\ref{equt:103}}^p \sup_{0\leq s\leq t}\sup_{x\in \IR^d}
\left(1+G_{\frac{\lambda}{2},p}u(s,x))\right).
\end{eqnarray}
Recall that
$a< \frac{1-\eta}{2}$
 and $p\geq p^{*}>\frac{1+\frac{d}{2}}{a}$. A bit of algebra shows that the 
whole expression in~(\ref{equt:115}) is finite and bounded by 
\[ C(T,\lambda,p)c_{\ref{equt:103}}^p \sup_{0\leq s\leq T}\sup_{x\in\IR^d}
\left(1+G_{\frac{\lambda}{2},p}u(s,x))\right).\]
This  together with~(\ref{equt:118}),~(\ref{equt:119})
proves~(\ref{equt:121}) for all $p\geq p^{*}$ with $\tilde{p}=p$.
Note, however, that
if $p< p^{*}$ then (\ref{equt:121}) also holds with $\tilde{p}=p^{*}$ due to the fact
that $u^{p} \leq 1+u^{p^{*}}$ for any $u\geq 0, p<p^{*}$. Hence we are
done. 

\end{proof}

The next result gives bounds on spatial and temporal differences of stochastic
convolution integrals which in particular will imply they are H\"older
continuous. The result is an adaptation of Theorem 2.1 of Sanz-Sol\'e and
Sarr\`a
\cite{SS02} to our situation.
\begin{lemma}
\label{thm:SanzSoleSarra}
Let $u$ be a solution to
(\ref{heatecolint}) satisfying the assumptions of Lemma \ref{lem:1}.
Define
\begin{equation}
\label{eq:stochintmn}
Z(t,x)= \int_0^t \int p_{t-s}(x-y)
\sigma(u(s,y)) W(ds dy),\; t\geq 0, \; x\in \IR^d.
\end{equation}
Then, for $T, R>0,$ and $0 \leq t,t'\leq T, x,x' \in \IR^{d}$ such that
  $|x-x'|<R$ as well as $p \in [2,\infty)$ and $\xi \in (0, 1-\eta)$
\begin{eqnarray}
\label{eq:Hoelderest}
\lefteqn{\IE\left( |Z(t,x)- Z(t',x')|^p e^{-\lambda |x|} \right)}
\\
\nonumber
&\leq& C(T,\lambda,p)c_{\ref{equt:103}}^p \left(1+\sup_{0\leq s\leq T}\sup_{z\in
\IR^d}  G_{\frac{\lambda}{p+1},p}u(s,z)\right)
 \left(|t-t'|^{\frac{\xi}{2} p} + |x-x'|^{\xi p}  \right).
\end{eqnarray}
In particular, if $G_{\frac{\lambda}{p+1},p}u(\cdot,\cdot)$ is bounded on 
$[0,T]\times\IR^d$, then there is a version of 
$Z$ which is uniformly
H\"older continuous on compact subsets of $[0,T]\times\IR^d$  with coefficients
$\frac{\xi}{2}$ in time and
$\xi$ in space.
\end{lemma}
\begin{proof}
The proof follows the proof of Theorem 2.1 in
Sanz-Sol{\'e} and Sarr{\`a} \cite{SS02}.
We use the same notation as in the
proof of Lemma \ref{lem:1}, so $Z(t,x) = J^{a-1}J_{a} u(t,x)$ by
(\ref{equt:111}). Now assume that $t'\geq t.$
By Lemma~\ref{lem:2} and H\"older's inequality we obtain 
\begin{eqnarray*}
& &\IE\left(|Z(t',x')- Z(t,x)|^p\right)e^{-\lambda |x|}\\
&\leq& C(p) \IE\Big( \Big| \int_{0}^{t} \int_{\IR^{d}}
(p_{t'-s}(x'-y)(t'-s)^{a-1} - p_{t-s}(x-y)(t-s)^{a-1}) \\
&&\phantom{AAAAAAAAAAAAAAAAAAAA}J_{a}u(s,y)
dy ds\Big|^{p} \Big)e^{-\lambda |x|}\\
&&+ C(p)\IE\left( \left|\int_{t}^{t'} \int_{\IR^{d}}
p_{t'-s}(x'-y)(t'-s)^{a-1}  J_{a}u(s,y) dy
ds\right|^{p}\right)e^{-\lambda |x|}\\
&\leq&  C(T,\lambda,p)c_{\ref{equt:103}}^p \left(1+\sup_{0\leq s\leq
T}\sup_{x\in
\IR^d} G_{\frac{\lambda}{p+1},p}u(s,x)\right)
 e^{-\lambda |x|}
\\
&&\mbox{}\times\Big\{
 \Big(\int_{0}^{t} \int_{\IR^{d}}
|p_{t'-s}(x'-y)(t'-s)^{a-1} - p_{t-s}(x-y)(t-s)^{a-1}| e^{\frac{\lambda}{p+1}|y|} dy
ds\Big)^{p} \\
&&\mbox{}+  \Big( \int_{t}^{t'} \int_{\IR^{d}}
p_{t'-s}(x'-y)(t'-s)^{a-1} e^{\frac{\lambda}{p+1}|y|} dy ds\Big) ^{p} \Big\}.
\end{eqnarray*}
Here, we have in the second inequality
also inserted additional factors of $e^{-\frac{\lambda}{p+1} |y|} e^{\frac{\lambda}{p+1}
|y|}$ so that we could apply Lemma~\ref{lem:2} to bound the expectation of
$J_{a}u$ by using Jensen's inequality.  From this point we proceed as in
\cite{SS02} (proof of Theorem 2.1), the only difference being that we have to
take care of the additional nuisance factors $e^{\frac{\lambda}{p+1} |y|}$.
This can be done with the help of (\ref{equt:104}) and (\ref{eq:pdiffest}) of
Lemma \ref{lemma:pdiffest} using the remaining factor $e^{-\lambda |x|}.$
\end{proof}

\vspace{1cm}

The next lemma assures that for any $u\in  C(\IR_{+}, C_{tem})$  which solves 
(\ref{heatecolint}), 
 $G_{\lambda,p}u(t,x)$ is bounded.

\begin{lemma}
\label{lem:5}
Let $u_{0}\in C_{tem}$ and  let $\sigma$ be a
continuous function satisfying the growth condition~(\ref{equt:103}).
 Assume that  $(A)_\eta$ holds for some $\eta \in [0,1).$  
If  $u\in  C(\IR_{+}, C_{tem})$ a.s.  is a solution to~(\ref{heatecolint}) then 
 it  satisfies
the following moment bound. For any $T>0$ and $p\geq 1$,
\begin{equation}
\label{equt:110}
\sup_{0 \leq t \leq T} \sup_{x \in \IR^{d}}
\IE \left( |u(t,x)|^{p} e^{-\lambda |x|} \right) \leq
C_{T,\lambda,p}(c_{\ref{equt:103}},\left\|u_{0}
 \right\|_{\frac{\lambda}{p},\infty}),
\end{equation}
where $C_{T,\lambda,p}(\cdot,\cdot)$ is bounded on the compacts of
 $\IR_+\times\IR_+$. 
\end{lemma}

\begin{proof}
Define 
\[ \tau_n= \inf\{t:\; \left\| u_{t}\right\|_{\frac{\lambda}{p},\infty}\geq
n\}  \]
 We set
\begin{equation*}
G^{\tau_{n}}_{\lambda,p}u(t,x):= \IE \left(|u(t,x)|^p 1(t\leq\tau_n)
e^{-\lambda |x|} \right).
\end{equation*}
Note that by definition,
\begin{eqnarray}
\label{equt:107} 
\sup_{s\leq t}\sup_{x\in\IR^d} G^{\tau_{n}}_{\lambda,p}u(t,x) \leq
n^p,\;\;\forall t\geq 0, n\geq 1. 
\end{eqnarray}
From~(\ref{heatecolint}) we get 
\begin{eqnarray}
\nonumber
&&\sup_{0\leq t\leq T}\sup_{x\in\IR^d}G^{\tau_{n}}_{\lambda,p}u(t,x) 
\\
&&\quad\leq
\label{eq:101}
C \sup_{0\leq t\leq T}\sup_{x\in\IR^d}\IE \left( 
\left|\int_{\IR^d} p_t(x-y) u_0(y) dy\right|^p e^{-\lambda |x|} \right)\\
\label{eq:102}
&&\qquad+ C\sup_{0\leq t\leq T}\sup_{x\in\IR^d}\IE \left(
|\int_0^{t} \int_{\IR^d} p_{t-s}(x-y) 
 1(s\leq \tau_n)\sigma(u(s,y)) W(dyds)|^p
e^{-\lambda |x|} \right).
\end{eqnarray}
By (\ref{equt:119}) the term on  line (\ref{eq:101}) is bounded by
\begin{eqnarray}
\label{equt:105}
 C(T,\lambda,p)\left\|u_{0}
 \right\|^p_{\frac{\lambda}{p},\infty}.
\end{eqnarray}

Again, as in Lemma~\ref{lem:1}, we use the factorization
method to bound the term in (\ref{eq:102}).
First, we assume that $p> \frac{2}{1-\eta} >2$ so that we can choose
a constant $a$ with $0< \frac{1}{p}<a<\frac{1-\eta}{2}<1.$
Recall~(\ref{equt:111}), and by 
 several applications of H\"older's inequality
we obtain 
\begin{eqnarray}
\nonumber
&&\IE \left( 
|J^{a-1}J_{a}^{\tau_n} u(t,x)|^{p} e^{-\lambda |x|}  \right)\\
\nonumber
&=& C\IE \left(
\left|\int_{0}^{t} \int_{\IR^{d}} (t-s)^{a-1}
p_{t-s}(x-y)J_{a}^{\tau_n}u(s,y) dy ds\right|^{p}
e^{-\lambda |x|}   \right)\\
\nonumber
&\leq& C\IE \left(
\left(\int_{0}^{t} (t-s)^{a-1} \left( \left|\int_{\IR^{d}}
p_{t-s}(x-y)J_{a}^{\tau_n}u(s,y) dy\right|^{p}  e^{-\lambda |x|}
\right)^{\frac{1}{p}}
\ ds \right)^{p}\right) \\
\nonumber
&\leq& C(T,\lambda) \IE \left(
\left(\int_{0}^{t} (t-s)^{a-1}
 \left(  \int_{\IR^{d}} |J_{a}^{\tau_n}u(s,y)|^{p}p_{t-s}(x-y)
e^{-\lambda |x|} dy \right)^{\frac{1}{p}}  \ ds \right)^{p} \right)\\
\nonumber
&\leq& C(T,\lambda)\left(\int_{0}^{T} s^{\frac{p}{p-1}(a-1)} ds \right)^{p-1}
 \cdot \left(\int_{0}^{t}  \int_{\IR^{d}} \IE \left( |J_{a}^{\tau_n}u(s,y)
|^{p}\right) e^{-\lambda |x|}p_{t-s}(x-y) dy ds\right) \\
\label{equt:106}
&\leq& C(T,\lambda,p)c_{\ref{equt:103}}^p
\left(1+ \int_{0}^{t} \sup_{0\leq r\leq s}\sup_{z\in \IR^d} G^{\tau_{n}}_{\lambda,p}u(r,z) ds\right),
\;\;\;\forall t\leq T, \;x\in\IR^d, 
\end{eqnarray}
where we have also used  Lemma~\ref{lem:2} and (\ref{equt:104}) in the last inequality
as well as $a>\frac{1}{p}.$
Taking (\ref{equt:105}) together with (\ref{equt:106}),
we obtain that  there is a constant
$C=C(T,\lambda,p)$ independent of $n$ such that for all $t \leq T,$
\begin{equation}
\sup_{0\leq s\leq t}\sup_{x\in\IR^d}G^{\tau_{n}}_{\lambda,p}u(t,x)
\leq C ( c_{\ref{equt:103}}^p+\left\|u_{0}
\right\|^p_{\frac{\lambda}{p},\infty})\left(
1+ \int_{0}^{t} \sup_{0\leq r\leq s}\sup_{x\in\IR^d}G^{\tau_{n}}_{\lambda,p}u(r,x)ds\right),\;\;\forall n\geq 1. 
\end{equation}
But the left hand side is bounded (due to (\ref{equt:107})).
Thus, by Gronwall's Lemma,
\begin{equation}
\label{equt:109}
\sup_{0\leq t\leq T}\sup_{x\in\IR^d}G^{\tau_{n}}_{\lambda,p}u(t,x)
 \leq C_{T,\lambda,p}(c_{\ref{equt:103}}, \left\|u_{0}
 \right\|_{\frac{\lambda}{p},\infty}), \forall n\geq 1, 
\end{equation}
where $C_{T,\lambda,p}(\cdot,\cdot)$ is bounded on the compacts of
 $\IR_+\times\IR_+$. 
(We have obtained this result with the restriction
$p>\frac{2}{1-\eta},$ which then
immediately implies that it is true for all $p>0$ since we are considering
$L^{p}$ norms with respect to a finite measure.)

Now, recall that $u\in  C(\IR_{+}, C_{tem})$ a.s. so that $\tau_n\uparrow\infty$ a.s., and
hence
\begin{eqnarray}
\nonumber
\IE \left(|u(t,x)|^p 
e^{-\lambda |x|}\right)
&=& 
\IE \left( \lim_{n\rightarrow \infty}|u(t,x)|^p 1(t\leq\tau_n)
e^{-\lambda |x|} \right)\\
\nonumber
&\leq& \liminf_{n\rightarrow \infty} G^{\tau_{n}}_{\lambda,p}u(t,x),
\end{eqnarray}
where the second inequality follows by Fatou's lemma. Use this and the fact
 that the right hand side of~(\ref{equt:109}) does not depend on $n$ to obtain
\begin{eqnarray}
\label{equt:116}
\sup_{0\leq t\leq T}\sup_{x\in\IR^d}\IE \left(|u(t,x)|^p 
e^{-\lambda |x|}\right)
&\leq& 
C_{T,\lambda,p}(c_{\ref{equt:103}}, \left\|u_{0}
\right\|_{\frac{\lambda}{p},\infty}), \forall n\geq 1, 
\end{eqnarray}
where $C_{T,\lambda,p}(\cdot,\cdot)$ is bounded on the compacts of
 $\IR_+\times\IR_+$. 

\end{proof}

\vspace{1cm}

\noindent
PROOF OF THEOREM \ref{thm:nonLipschitzexistence}.
Recall our hypotheses imply $(A)_\eta$  holds for some $\eta\in[0,1)$ (see
Remark~\ref{rem:sm}). We can choose a sequence of Lipschitz continuous functions
$\sigma_{n}$ on
$\IR^{d}$ such that the growth bound (\ref{growthcond}) holds uniformly
($\sigma_{n}(u) \leq c_{\ref{growthcond}}(1+ |u|)$ for all $u \in \IR, n \in
\IN$), and such that the
$\sigma_{n}$ converge uniformly to $\sigma$ as $n \rightarrow \infty.$
We also set
\begin{equation*}
u^{m}_{0}(x)=
\left\{\begin{array}{ll}
u_{0}(x) & \text{ if } |u_{0}(x)|< m,\\
m & \text{ if } u_{0}(x)\geq  m,\\
-m & \text{ if } u_{0}(x)\leq  -m,
\end{array}\right.
\end{equation*}
which implies that $u^{m}_{0}\in C_{b}(\IR^{d})$ and 
\begin{equation}
\label{eq:minitcond}
\sup_{m \in \IN} \sup_{x \in \IR^{d}} |u_0^m(x)| e^{-\lambda |x|} < \infty. 
\end{equation}
Hence, by Theorem \ref{thm:Dalang}, for each $m,n$ there exists a unique
solution to 
\begin{eqnarray}
\label{equt:123}
u^{m,n}(t,x) &=& \int_{\IR^{d}} p_{t}(x-y) u^{m}_{0}(y) dy
+ \int_{0}^{t}  \int_{\IR^{d}} p_{t-s}(x-y) \sigma_{n}(u^{m,n}(s,y)) W(dy ds).
\end{eqnarray}
It is easy to check that the first term on the right hand side of~(\ref{equt:123}) is jointly continuous
 on $[0,\infty) \times \IR^{d}$.
Moreover, by Theorem \ref{thm:Dalang}, 
$\sup_{0 \leq t \leq T} \sup_{x \in \IR^{d}}
\IE\left( |u^{m,n}(t,x)|^{p} \right) < \infty$, and hence by Lemmas~\ref{lem:1},~\ref{thm:SanzSoleSarra}
we obtain that  
\begin{eqnarray}
 u^{m,n}\in C(\IR_{+}\,, C_{tem}).
\end{eqnarray}
Now let us go to the limit as $m,n\rightarrow \infty$. Let $Z^{m,n}$ denote that stochastic integral 
 on the right hand side of~(\ref{equt:123}). Since $u^{m,n}\in C(\IR_{+}\,, C_{tem})$ we may apply 
 Lemmas~\ref{thm:SanzSoleSarra},~\ref{lem:5} to get 
\begin{eqnarray}
\label{equt:124}
\IE\left( |Z^{m,n}(t,x)- Z^{m,n}(t',x')|^p e^{-\lambda |x|} \right)
&\leq& C_{T,\lambda,p}(c_{\ref{growthcond}},\left\|u^{m}_{0}
\right\|_{\frac{\lambda}{p},\infty}),
\end{eqnarray}
where $C_{T,\lambda,p}(\cdot,\cdot)$ is bounded on the compacts of
 $\IR_+\times\IR_+$. (\ref{equt:124}) combined with a Kolmogorov type
tightness criterion (see Lemma 6.3 of Shiga \cite{tS94}) now implies
that the stochastic integrals $Z^{m,n}$ 
are tight in $C(\IR_{+}, C_{tem}).$ It is also not hard to show that
$(t,x)\mapsto \int_{\IR^{d}} p_{t}(x-y) u_{0}^{m} dy$ are tight in $C(\IR_{+},
C_{tem})$ by using the Arzela-Ascoli Theorem and the uniformity in $m$ as in
(\ref{eq:minitcond}). 

Therefore, $u^{m,n}$ are tight in $C(\IR_{+}, C_{tem})$ and we can choose
an appropriate probability space and define $u^{m,n}$ on it identical in distribution
to a subsequence of the
original sequence of solutions which converge a.s. in $C(\IR_{+}, C_{tem})$ to
some process $u.$
It is routine to establish from this that all the terms in (\ref{heatecolint})
converge a.s. to the appropriate limits so that the limit $u\in C(\IR_{+},
C_{tem})$ is indeed a solution to (\ref{heatecolint}) with the desired $\sigma$
and $u_{0}
\in C_{tem}.$

\hfill$\mathbf{\Box}$ \vspace{.5\baselineskip}

\vspace{0.5cm}

\noindent
PROOF OF PROPOSITION \ref{prop}.
Recall from Remark~\ref{rem:sm} that our hypotheses imply $(A)_\eta$ for any
$\eta\in (\alpha/2,1)$.  Lemmas~\ref{lem:1} and \ref{lem:5} now imply (a). 
Now use Lemma~\ref{thm:SanzSoleSarra} (and
Lemma~\ref{lem:5}) to derive part (b) for $Z$.
For $u_0\in C_{tem}$, $S_tu_0(x)\equiv\int p_t(y-x)u_0(y)\,dy$ is smooth on
$(0,\infty)\times\IR^d$, and so is uniformly Lipschitz on compact subsets of
$(0,\infty)\times\IR^d$.  This gives the required H\"older continuity for $u$.  
\hfill$\mathbf{\Box}$ \vspace{.5\baselineskip}

{\small
\bibliographystyle{alpha}
\bibliography{references}}

\end{document}